\documentclass[letter]{article}

\usepackage[numbers]{natbib}
\usepackage{fullpage}
\usepackage{amsmath,amssymb,amsthm}
\usepackage{tikz}
\usepackage{pgfplots}
\usepackage{booktabs}
\usepackage{units}
\usepackage{subfig}

\usepackage{algorithm}
\usepackage{algpseudocode}

\usepackage{multirow}

\newcommand{\Mat}[1]{\mathsf{#1}}  
\newcommand{\op}[1]{\mathcal{#1}}  

\theoremstyle{plain}

\graphicspath{{pics/}}

\begin{document}

\title{Optimized M2L Kernels for the Chebyshev Interpolation based Fast
  Multipole Method} \author{Matthias Messner, Berenger Bramas, Olivier
  Coulaud, Eric Darve} \date{\today}
\maketitle

\begin{abstract}
A fast multipole method (FMM) for asymptotically smooth kernel functions
($\nicefrac{1}{r}$, $\nicefrac{1}{r^4}$, Gauss and Stokes kernels, radial
basis functions, etc.) based on a Chebyshev interpolation scheme has been
introduced in \citep{fong09a}. The method has been extended to oscillatory
kernels (e.g., Helmholtz kernel) in \citep{Messner2011}. Beside its generality
this FMM turns out to be favorable due to its easy implementation and its high
performance based on intensive use of highly optimized BLAS
libraries. However, one of its bottlenecks is the precomputation of the
multiple-to-local (M2L) operator, and its higher number of floating point
operations (flops) compared to other FMM formulations. Here, we present
several optimizations for that operator, which is known to be the costliest
FMM operator. The most efficient ones do not only reduce the precomputation
time by a factor up to $340$ but they also speed up the matrix-vector
product. We conclude with comparisons and numerical validations of all
presented optimizations.


\end{abstract}


\section{Introduction}
\label{sec:introduction}

The fast multipole method (FMM) is a method first designed in
\citep{greengardRokhlin:87} to reduce the cost of matrix-vector products from
$\op{O}(N^2)$ to $\op{O}(N)$ or $\op{O}(N \log N)$ depending on the underlying
kernel function. Most FMM variants have been developed and optimized for
specific kernel functions \citep{Rokhlin93, greengardRokhlin:97, DarveEtAl04,
  ChengEtAl06}. However, some have a formulation that is independent of the
kernel function \citep{Ying2004591, Martinsson:2007, Gimbutas:2002, fong09a,
  Messner2011}. The paper at hand addresses the optimization of one of these
formulations, the so called black-box FMM ({\bf bbFMM}), presented in
\citep{fong09a}. It is based on the approximation of the kernel function via a
Chebyshev interpolation and is a black-box scheme for kernel functions that
are asymptotically smooth, e.g., $1/(r^2+c^2)^{n/2}$ with $r=|x-y|$, $c$ a
given parameter and $n\in\mathbb{N}$. The bbFMM has been extended to the
directional FMM ({\bf dFMM}) for oscillatory kernels in
\citep{Messner2011}. It is suitable for any kernel function of the type $g(r)
e^{\imath k r}$ where $g(r)$ is an asymptotically smooth function ($\imath^2 =
-1$ is the imaginary unit and $k$ the wavenumber).

The main idea of the FMM is to properly separate near-field $(|x-y|
\rightarrow 0)$ and far-field $(|x-y| \rightarrow \infty)$. The near-field is
evaluated directly and the far-field can be approximated and thus computed
efficiently. In this paper, we will denote M2L the operator that transforms a
multipole expansion into a local expansion. The M2L operator is the costliest
step in the method, since it needs to be applied to all clusters in the
interaction list, that is $189$ times for each cluster for a typical
configuration in bbFMM. In this paper we focus on various
optimizations of this operator for both the bbFMM and the dFMM.

First, we address the optimization proposed in \citep{fong09a,
  Messner2011}. In that paper, the M2L operator between a cluster and its
interaction list is further compressed using a singular value decomposition
(SVD). The key idea is that the SVD provides the smallest possible rank given
an error $\varepsilon$ (therefore leading to the smallest computational
cost). In \citep{fong09a}, a single SVD is used to compress {\bf all the M2L
  operators at a given level} using a single singular vector basis. However,
the singular values of individual M2L operators decay at different rates; it
is often the case for example that cluster pairs that are separated by a small
distance have higher rank than pairs that are separated by a greater
distance. If we denote $w$ the diameter (or width) of a cluster and $d$ the
distance between two clusters, then the smallest distance corresponds to
$d_\text{min} = 2w$ while the largest is $d_\text{max} = 3\sqrt{3} w$. The
ratio $d_\text{max} / d_\text{min} = 3\sqrt{3}/2$ is in fact relatively
large. This can be taken advantage of to further compress the M2L operator on
a {\bf cluster-pair basis} leading to an even smaller computational cost than
the original algorithm \citep{fong09a}. This point is investigated in details
in this paper.

Another bottleneck in the FMM of \citep{fong09a} is the precomputation time of
the M2L operators. We therefore introduce a new set of optimizations that
exploit symmetries in the arrangement of the M2L operators. For example, for
bbFMM, this allows us to express all M2L operators ($316 = 7^3-3^3$) as
permutations of a subset with only $16$ unique M2L operators. These
permutations correspond to reflections across various planes (for example
$x=0$, $y=0$, $z=0$, etc). Table~\ref{tab:prectimings} for example reports
reductions in precomputation time by a factor of about $50$--$200$x.

Besides drastically reducing precomputation time and memory requirement this
approach paves also the road for various algorithmic improvements. Modern
processors use memory cache to enhance performance and, as a result, schemes
that block data (thereby reducing the amount of data traffic between main
memory and the floating point units) can result in much higher performance. In
our case, the use of symmetries allows calling highly optimized matrix-matrix
product implementations (instead of matrix-vector) that have better cache
reuse. Overall, this may result in acceleration by a factor $4$--$6$x for
smooth kernels (Tables~\ref{tab:timebbfmm_e0}--\ref{tab:timebbfmm_e2}). We
also present results using oscillatory kernels of the type mentioned above
($g(r) e^{\imath k r}$). In that case, the acceleration is more modest and
ranges from about $1.2$ to $2.7$x
(Table~\ref{tab:timem2l_e0}--\ref{tab:timem2l_e2}).

In this paper we therefore focus both on performance and precomputation
time. Both are important factors for the FMM: the precomputation is crucial if
we are interested in only one matrix-vector product, while in other cases,
such as the solution of a linear system, the fast application of matrix-vector
products is central.

The paper is organized as follows. In Sec.~\ref{sec:fast_summation_scheme}, we
briefly recall the bbFMM and the dFMM and introduce notations that are needed
for explanations later in this paper. In Sec.~\ref{sec:m2l_operators}, we
address the separation of near- and far-field, introduce the notion of
transfer vector to uniquely identify M2L operators, and describe how the
kernel (smooth or oscillatory) affects the interaction list. We start
Sec.~\ref{sec:m2l_optimizations} with a brief recall of the known optimization
of the M2L operator (see \citep{fong09a, Messner2011}) and suggest further
improvements. Then, we present a new set of optimizations and explain them in
details. These include using a low-rank approximation of M2L for individual
interactions (individual cluster pairs), along with the definition of
symmetries and permutations and how they are used to reduce the computational
cost. Finally, in Sec.~\ref{sec:results}, we present numerical results. We
compare all variants for bbFMM and dFMM with three datasets corresponding to a
sphere, an oblate spheroid, and a prolate spheroid. The measured running time
of the FMM, along with the precomputation time, and its accuracy as a function
of various parameters are reported. The efficiency of the algorithms as a
function of the target accuracy is considered. We also performed tests using
two different kinds of linear algebra libraries for the matrix-matrix products
and other operations (an implementation of BLAS and LAPACK vs the Intel MKL
library).


\section{Pairwise particle interactions}
\label{sec:fast_summation_scheme}

The problem statement reads as follows. Assume, the cells $X$ and $Y$ contain
source $\{y_j\}_{j=1}^N$ and target particles $\{x_i\}_{i=1}^M$,
respectively. Compute the interactions
\begin{equation}
  \label{direct}
  f_i = \sum_{j=1}^N K(x_i,y_j)w_j \quad \text{for } i,\dots,M.
\end{equation}
The kernel function $K(x,y)$ describes the influence of the source particles
onto the target particles. The cost of directly evaluating the summation in
Eqn.~\eqref{direct} grows like $\op{O}(MN)$ which becomes prohibitive as $M,N
\rightarrow \infty$ and it is why we need a fast summation scheme.

\subsection{Fast summation schemes based on Chebyshev interpolation}

For a detailed derivation and error analysis of the FMM based on Chebyshev
interpolation we refer the reader to \citep{fong09a, Messner2011}. We adapt
most of their notations and repeat only concepts which are necessary to
understand explanations hereafter.

Let the function $f : \mathbb{R}^3 \rightarrow \mathbb{C}$ be approximated by
a Chebyshev interpolation scheme as
\begin{equation}
  \label{eq:chebinter1d}
  f(x) \sim \sum_{|\alpha|\le 3\ell} S_\ell(x,\bar x_\alpha) f(\bar x_\alpha) 
\end{equation}
with the $3$-dimensional multi-index $\alpha = (\alpha_1, \alpha_2, \alpha_3)$
and $|\alpha| = \max(\alpha_1, \alpha_2, \alpha_3)$ with $\alpha_i \in
(1,\dots,\ell)$. The interpolation points $\bar x_\alpha = (\bar x_{\alpha_1},
\bar x_{\alpha_2}, \bar x_{\alpha_3})$ are given by the tensor-product of the
Chebyshev roots $\bar x_{\alpha_i}$ of the Chebyshev polynomial of first kind
$T_\ell(x) = \cos (\arccos x)$ with $x \in [-1,1]$. The interpolation operator
reads as
\begin{equation}
  \label{eq:interpolop3d}
  S_\ell(x,\bar x_\alpha) = S_\ell(x_1,\bar x_{\alpha_1}) \,
  S_\ell(x_2,\bar x_{\alpha_2}) \, S_\ell(x_3,\bar x_{\alpha_3}).
\end{equation}
For the interpolation on arbitrary intervals, we need the affine mapping $\Phi
: [-1,1] \rightarrow [a,b]$. We omit it hereafter for the sake of readability.

\subsubsection{Black--box FMM (bbFMM)}
If two cells $X$ and $Y$ are well separated, we know from \citep{fong09a} that
\emph{asymptotically smooth kernel functions} can be interpolated as
\begin{equation}
  \label{eq:interpolated_laplace}
  K(x,y) \sim \sum_{|\alpha|\le \ell} S_\ell(x,\bar x_\alpha)
  \sum_{|\beta|\le \ell} K(\bar x_\alpha, \bar y_\beta) \; S_\ell(y,\bar
  y_\beta).  
\end{equation}
We insert the above approximation into Eqn.~\eqref{direct} and obtain
\begin{equation}
  \label{eq:fast_summation_laplace}
  f_i \sim \sum_{|\alpha|\le \ell} S_\ell(x_i,\bar x_\alpha) \sum_{|\beta|\le \ell}
  K(\bar x_\alpha, \bar y_\beta) \sum_{j=1}^{N} S_\ell(y_j,\bar y_\beta) \; w_j
\end{equation}
which we split up in a three-stage fast summation scheme.
\begin{enumerate}
\item Particle to moment (P2M) or moment to moment (M2M) operator: equivalent
  source values are anterpolated at the interpolation points $\bar y_\beta \in
  Y$ by
  \begin{equation}
    \label{eq:m2m_bbfmm}
    W_\beta = \sum_{j=1}^N S(y_j, \bar y_\beta) \, w_j
    \quad \text{for} \quad |\beta|\le \ell. 
  \end{equation}
\item Moment to local operator (M2L): target values are evaluated at the
  interpolation points $\bar x_\alpha \in X$ by
  \begin{equation}
    \label{eq:m2l_bbfmm}
    F_\alpha = \sum_{|\beta|\le \ell} K(\bar x_\alpha,\bar y_\beta) \,
    W_\beta \quad \text{for} \quad |\alpha|\le \ell. 
  \end{equation}
\item Local to local (L2L) or local to particle (L2P) operator: target values
  are interpolated at final points $x_i\in X$ by
  \begin{equation}
    \label{eq:l2l_bbfmm}
    f_i \sim \sum_{|\alpha|\le \ell} S(x_i,\bar x_\alpha) \, F_\alpha
    \quad \text{for} \quad i=1,\dots,M.  
  \end{equation}
\end{enumerate}
Recall, the cells $X$ and $Y$ are well separated and thus all contributions of
$f_i$ can be computed via the above presented fast summation scheme (no direct
summation is necessary).

\subsubsection{Directional FMM (dFMM)}
Whenever we deal with \emph{oscillatory kernel functions}, such as the
Helmholtz kernel, the wavenumber $k$ comes into play. Depending on the
diameter of the cells $X$ and $Y$ and the wavenumber they are either in the
low-frequency or in the high-frequency regime. In the low-frequency regime the
fast summation schemes of the dFMM and the bbFMM are the same. In the
high-frequency regime the fast summation scheme becomes directional.
From \citep{Messner2011} we know that any oscillatory kernel function $K(x,y)
= G(x,y) e^{\imath k |x-y|}$, where $G(x,y)$ is an asymptotically smooth
function, can be rewritten as
\begin{equation}
  \label{eq:helmkerneldir}
  K(x,y) =  K^u(x,y) e^{\imath k u \cdot (x-y)}
  \quad \text{with} \quad
  K^u(x,y) = G(x,y) e^{\imath k (|x-y| - u \cdot (x-y))}.
\end{equation}
We assume that the cells $X$ and $Y$ of width $w$ are centered at $c_x$ and
$c_y$ and $c_y$ lies in a cone of direction $u$ being centered at $c_x$ (think
of the domain around $X$ being virtually subdivided in cones given by
directional unit vectors $\{u_c\}_{c=1}^C$, where $C$ is determined by their
aperture). If the cell pair $(X,Y)$ satisfies the \emph{separation} criterion
$\op{O}(kw^2)$ and the \emph{cone-aperture} criterion $\op{O}(1/kw)$, the
error of the Chebyshev interpolation of the kernel function
\begin{equation}
  \label{eq:interp_kernel}
  K^u(x,y) \sim \sum_{|\alpha|\le \ell} S_\ell(x,\bar x_\alpha)
  \sum_{|\beta|\le \ell} K^u(\bar x_\alpha, \bar y_\beta) \, S_\ell(y,\bar
  y_\beta)
\end{equation}
decays exponentially in the interpolation order $\ell$ (independent of the
wavenumber $k$; see~\citep{mason2003chebyshev, Messner2011}). We insert the
above interpolated kernel function in Eqn.~\eqref{direct} and obtain
\begin{equation}
  \label{eq:fast_sum}
  f_i \sim e^{\imath k u\cdot x_i} \sum_{|\alpha|\le \ell} S_\ell(x_i,\bar
  x_\alpha) e^{-\imath k u\cdot \bar{x}_\alpha} \sum_{|\beta|\le \ell}
  K(\bar x_\alpha, \bar y_\beta) e^{\imath k u\cdot \bar{y}_\beta}
  \sum_{j=1}^{N} S_\ell(y_j,\bar y_\beta) e^{-\imath k u\cdot y_j}\, w_j 
\end{equation}
for all $i=1,\dots,M$. Similarly as with the bbFMM, a three-stage fast
summation scheme for oscillatory kernels in the high-frequency regime can be
constructed.
\begin{enumerate}
\item Particle to moment (P2M) or moment to moment (M2M) operator: equivalent
  source values are anterpolated at the interpolation points $\bar y_\beta \in
  Y$ by
  \begin{equation}
    \label{eq:m2m_dfmm}
    W_\beta^u = e^{\imath k u\cdot \bar{y}_\beta} \sum_{j=1}^N S(y_j, \bar y_\beta) \,
    e^{-\imath k u\cdot y_j} \, w_j \quad \text{for} \quad |\beta|\le \ell. 
  \end{equation}
\item Moment to local operator (M2L): target values are evaluated at the
  interpolation points $\bar x_\alpha \in X$ by
  \begin{equation}
    \label{eq:m2l_dfmm}
    F_\alpha^u = \sum_{|\beta|\le \ell} K(\bar x_\alpha,\bar y_\beta) \,
    W_\beta^u \quad \text{for} \quad |\alpha|\le \ell. 
  \end{equation}
\item Local to local (L2L) or local to particle (L2P) operator: target values
  are interpolated at final points $x_i\in X$ by
  \begin{equation}
    \label{eq:l2l_dfmm}
    f_i \sim e^{\imath k u\cdot x_i} \sum_{|\alpha|\le \ell} S(x_i,\bar
    x_\alpha) \, e^{-\imath k u\cdot \bar{x}_\alpha} \, F_\alpha^u \quad
    \text{for} \quad i=1,\dots,M.  
  \end{equation}
\end{enumerate}

Even though the bbFMM and the dFMM are here presented as single-level schemes,
they are usually implemented as multilevel schemes. Strictly speaking, the
steps one and three of both schemes are the P2M and L2P operators. Let us
recall briefly on the directional M2M and L2L operators of dFMM: based on the
criterion $\op{O}(1/kw)$, the aperture of the cones at the upper level is
about half the aperture at the lower level. Due to a nested cone construction
along octree levels, we preserve the accuracy of the Chebyshev interpolation
within the multilevel scheme. For a detailed description of all operators we
refer to \citep{fong09a, Messner2011}). Note, the similarity of the M2L
operators (step two of both schemes) of bbFMM and dFMM. In fact, the only
difference in their implementation is that in the bbFMM case we have one loop
over all cell pairs, whereas in the dFMM case we have two loops: the outer
loop over all existing cones of direction $\{u_c\}_{c=1}^C$ and the inner loop
over all cell pairs lying in the current cone. In this paper, we focus on the
M2L operators and their efficient numerical treatment.


\section{M2L operators}
\label{sec:m2l_operators}

The first step of any FMM consists in a proper separation of near- and
far-field. After that, the near-field is evaluated directly and the far-field
efficiently using a fast summation scheme. In this section, we focus on the
first step. The union of near- and far-field of a target cell $X$ is spatially
restricted to the near-field of its parent
cell. Algorithm~\ref{alg:near_far_field} explains how these interactions are
computed for dFMM \citep{Messner2011}. The recursive partitioning starts with
the two root cells $X$ and $Y$ of the octrees for source and target
particles. If a pair of cells satisfies the separation criterion in the high-
or low-frequency regime, $Y$ is a far-field interaction of $X$. Else, if they
are at the leaf level of the octree, $Y$ is a near-field interaction of
$X$. If none is true, the cell is subdivided and the tests are repeated
recursively.
\begin{algorithm}
  \caption{Separate near- and far-field in the low- and high-frequency regime}
  \label{alg:near_far_field}
  \begin{algorithmic}[1]
    \Function{SeparateNearAndFarField}{$X,Y$} 
    \If{$(X,Y)$ are admissible in the high-frequency regime} 
    \State{add $Y$ to the directional far-field of
      $X$} \label{alg:ln:directional} 
    \Return
    \ElsIf{$(X,Y)$ are admissible in the low-frequency regime} 
    \State{add $Y$ to the far-field of $X$}
    \Return
    \ElsIf{$(X,Y)$ are leaves}
    \State{add $Y$ to the near-field of $X$}
    \Return
    \Else
    \For{all $X_{\text{child}} \in X$ and all $Y_{\text{child}} \in Y$}
    \State \Call{SeparateNearAndFarField}{$X_{\text{child}},
      Y_{\text{child}}$} 
    \EndFor
    \EndIf
    \EndFunction
  \end{algorithmic}
\end{algorithm}
In line~\ref{alg:ln:directional} in Alg.~\ref{alg:near_far_field} we use the
term directional far-field, a concept explained in detail in
\citep{Messner2011}: In the high-frequency regime the far-field is subdivided
into cones of direction $u$ needed by the directional kernel function
$K^u(x,y)$. Each source cell $Y$ is assigned to a cone and there are as many
directional far-fields as there are cones.

\subsection{Transfer vectors}
\label{sec:transfer_vectors}

In order to address interactions uniquely we introduce transfer vectors $t =
(t_1,t_2,t_3)$ with $t\in\mathbb{Z}^3$. They describe the relative positioning
of two cells $X$ and $Y$ and are computed as $t = \nicefrac{(c_x-c_y)}{w}$
where $c_x$ and $c_y$ denote the centers of the cells $X$ and $Y$ and $w$
their width. In the following we use transfer vectors to uniquely identify the
M2L operator that computes the interaction between a target cell $X$ and a
source cell $Y$. In matrix notation an M2L operator reads as $\Mat{K}_t$ of
size $\ell^3 \times \ell^3$ and the entries are computed as
\begin{equation}
  \label{eq:m2lmatrix}
  (\Mat{K}_t)_{m(\alpha) n(\beta)} = K(\bar x_\alpha, \bar y_\beta)
\end{equation}
with the interpolation points $\bar x_\alpha \in X$ and $\bar y_\beta \in
Y$. Bijective mappings
\begin{equation}
  \label{eq:mn_bijective}
  m,n : \{1,\dots,\ell\} \times \{1,\dots,\ell\} \times \{1,\dots,\ell\}
  \rightarrow \{1,\dots,\ell^3\}
\end{equation}
with $m^{-1}(m(\alpha))=\alpha$ and $n^{-1}(n(\beta))=\beta$ provide unique
indices to map back and forth between multi-indices of $\bar x_\alpha$ and
$\bar y_\beta$ and rows and columns of $\Mat{K}_t$. We choose them to be
$m(\alpha) = (\alpha_1-1) + (\alpha_2-1) \ell + (\alpha_3-1) \ell^2 + 1$ and
same for $n(\beta)$.

\subsection{Interaction lists}
\label{sec:translation_invariance}

Very commonly fast multipole methods are used for translation invariant kernel
functions $K(x,y) = K(x+v,y+v)$ for any $v\in\mathbb{R}^3$. Because of that
and because of the regular arrangement of interpolation points $\bar x$ and
$\bar y$ in uniform octrees it is sufficient to identify unique transfer
vectors at each level of the octree and to compute the respective M2L
operators. In the following we refer to such sets of unique transfer vectors
as interaction lists $T \subset \mathbb{Z}^3$.

If we consider \emph{asymptotically smooth kernel functions} the near-field is
limited to transfer vectors satisfying $|t| \le \sqrt{3}$; it leads to $3^3 =
27$ near-field interactions (see \citep{fong09a}). In a multi-level scheme,
these $27$ near-field interactions contain all $6^3=216$ near- and far-field
interactions of its eight child-cells. Far-field interactions are given by
transfer-vectors that satisfy $|t| > \sqrt{3}$. This leads to a maximum of
$6^3-3^2 = 189$ interactions per cell and we end up with the usual overall
complexity of $\op{O}(N)$ of fast multipole methods for asymptotically smooth
kernel functions. The union of all possible far-field interactions of the
eight child cells gives $7^3-3^3 = 316$ interactions. That is also the largest
possible number of different M2L operators have to be computed per octree
level. Most asymptotically smooth kernel functions happen also to be
homogeneous $K(\alpha r) = \alpha^n K(r)$ of degree $n$. In other words, if we
scale the distance $r=|x-y|$ between source and target by a factor of $\alpha$
the resulting potential is scaled by $\alpha^n$, where $n$ is a constant that
depends on the kernel function. The advantage of homogeneous kernel functions
is that the M2L operators need to be precomputed only at one level and can
simply be scaled and used on other levels. This affects the precomputation
time and the required memory.

If we consider \emph{oscillatory kernel functions} we need to distinguish
between the low- and high-frequency regime \citep{Messner2011}. The
admissibility criteria in the low-frequency regime are the same as those for
asymptotically smooth kernel functions. However, in the high-frequency regime
the threshold distance between near-field and far-field is $\op{O}(kw^2)$;
nonetheless, as shown in \citep{Messner2011} we end up with the usual
complexity of $\op{O}(N \log N)$ of fast multipole methods for oscillatory
kernel functions. It depends on the wavenumber $k$ (a property of the kernel
function). Thus, the size of near- and far-field is not a priori known as it
is in the case of asymptotically smooth kernel
functions. Table~\ref{tab:near_far_field_size} summarizes the number of near
and far-field interactions to be computed depending on different kernel
functions.
\begin{table}[htbp]
  \caption{Number of near- and far-field interactions depends on the kernel
    function}  
  \label{tab:near_far_field_size}
  \centering
    \begin{tabular}{l | c c }
      \toprule
      type of kernel function & cells in near-field & cells in far-field \\ 
      \midrule
      smooth                 & $\le 27$ per leaf  & $\le 316$ per level \\
      smooth and homogeneous & $\le 27$ per leaf  & $\le 316$ for all levels \\
      oscillatory            & depends on $k$ & depends on $k$ \\
      \bottomrule
  \end{tabular}
\end{table}


\section{Optimizing the M2L operators}
\label{sec:m2l_optimizations}

In all fast multipole methods the M2L operator adds the largest contribution
to the overall computational cost: for bbFMM it grows like $\op{O}(N)$ and for
dFMM like $\op{O}(N \log N)$. In this section, we first briefly recall the
optimization that was used up to now and suggest an improvement. Then, we
present a new set of optimizations that exploit the symmetries in the
arrangement of the M2L operators.

\subsection{Single low-rank approximation (SA)}
\label{sec:fast_convolution_former}

In the following we explain the basics of the optimization used in
\citep{fong09a, Messner2011} and we refer to it as the SA variant
hereafter. The idea is based on the fact that all M2L operators $\Mat{K}_t$
with $t \in T$ can be assembled as a single big matrix in two ways: either as
a row of matrices $\Mat{K}^{(\text{row})} = [\Mat{K}_1, \dots, \Mat{K}_t,
\dots, \Mat{K}_{|T|}]$ or as a column of matrices $\Mat{K}^{(\text{col})} =
[\Mat{K}_1; \dots; \Mat{K}_t; \dots; \Mat{K}_{|T|}]$ of M2L operators. The
cardinality $|T|$ gives the number of transfer vectors in the interaction list
$T$. Next, both big matrices are compressed using truncated singular value
decompositions (SVD) of accuracy $\varepsilon$ as
\begin{equation}
  \label{eq:2bigsvds}
  \Mat{K}^{(\text{row})} \sim \Mat{U\Sigma V}^* \quad \text{and} \quad
  \Mat{K}^{(\text{col})} \sim \Mat{A\Gamma B}^*
\end{equation}
with the unitary matrices $\Mat{U,B}$ of size $\ell \times r$ and $\Mat{V,A}$
of size $|T|\ell^3 \times r$ and the $r$ singular values in
$\Mat{\Sigma,\Gamma}$. With a few algebraic transformations each M2L operator
can be expressed as
\begin{equation}
  \label{eq:ucb_no_compression}
  \Mat{K}_t \sim \Mat{UC}_t \Mat{B}^*
  \quad \text{where} \quad
  \Mat{C}_t = \Mat{U}^*\Mat{K}_t\Mat{B} \quad \text{of size } r \times r
  \text{ is computed as} \quad
  \Mat{C}_t = \Mat{\Sigma V}_t^* \Mat{B}
  \quad \text{or} \quad
  \Mat{C}_t = \Mat{U}^* \Mat{A}_t \Mat{\Gamma}.
\end{equation}
The advantage of this representation is that the cost of applying the M2L
operator gets reduced from applying a matrix of size $\ell^3 \times \ell^3$ to
a matrix of only $r \times r$. Moreover, less memory is required. However, the
precomputation time grows cubically with the accuracy of the method due to the
complexity of the SVD. In \citep{Messner2011} the SVD has been substituted by
the adaptive cross approximation (ACA) followed by a truncated SVD
\citep{bebendorf:05}. The precomputation time has been cut down drastically
due to the linear complexity of the ACA.

\subsubsection{SA with recompression (SArcmp)}
If we use the SA variant the achieved low-rank $r$ is the same for all M2L
operators given by $\Mat{C}_t$. To a large extent, $r$ is determined by the
greatest individual low-rank of the M2L operators. This means that most of the
matrices $\Mat{C}_t$ of size $r \times r$ have effectively a smaller low-rank
$r_t \le r$. We exploit this fact by individually approximating them as
\begin{equation}
  \label{eq:recompression}
  \Mat{C}_t \sim \bar{\Mat{U}}_t \bar{\Mat{V}}_t^* \quad \text{with }
  \bar{\Mat{U}}_t,\bar{\Mat{V}}_t \text{ of size } r \times r_t
  \text{ and the constraint} \quad r_t < r/2.
\end{equation}
Without the constraint the low-rank representation is less efficient than the
original representation.
The effects of the recompression are studied in
Sec.~\ref{sec:recompression_ucb}.

\subsection{Individual low-rank approximation (IA)}
\label{sec:approach2}
As opposed to the SA approach, an individual low-rank approximation of the M2L
operators as
\begin{equation}
  \label{eq:ialowrank}
  \Mat{K}_t \sim \Mat{U}_t \Mat{V}_t^* \quad \text{with } \Mat{U}_t \Mat{V}_t
  \text{ of size } \ell^3 \times r_t
\end{equation}
directly leads to the optimal low-rank representation of each of them. As in
the previous section, the approximation can be performed by either a truncated
SVD or the ACA followed by a truncated SVD (note, the rank $r_t$ in the
Eqns.~\eqref{eq:recompression} and \eqref{eq:ialowrank} might be similar but
has not to be the same). All these variants (SA, SArcmp and IA) still require
the approximation and storage of all M2L operators. In terms of time and
memory however, it would be desirable to come up with a method that requires
the approximation and the storage of a subset of operators only. Let us
present a set of optimizations that fulfill these two requests.

\subsubsection{Symmetries and permutations}
\label{sec:symmetries}

Here, we illustrate how the full set of M2L operators can be expressed by a
subset only. The idea is based on symmetries in the arrangement of M2L
operators and exploits the uniform distribution of the interpolation
points. We start by presenting the idea using a model example. After
generalizing this idea, we demonstrate that M2L operators can be expressed as
permutations of others.

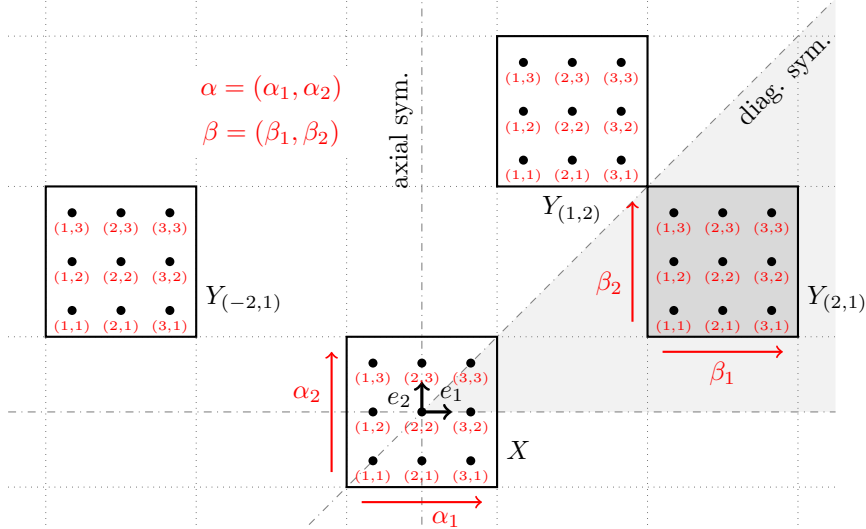
\begin{figure}[h!]
  \centering
  \begin{tikzpicture}
    \draw[very thin,dotted,step=2] (-6.5,-2.5) grid (4.5,4.5);

    \fill[gray, fill opacity=.1] (-1,-1) -- (4.5,-1) -- (4.5,4.5);

    \draw[gray,dashdotted] (-6.5,-1) -- (4.5,-1);
    \draw[gray,dashdotted] (-1,-2.5) to node[near end,sloped,above,black]
    {axial sym.} (-1,4.5);
    \draw[gray,dashdotted] (-2.5,-2.5) to node[very near end,sloped,below,black]
    {diag. sym.} (4.5,4.5);

    \draw[very thick,->] (-1,-1) -- (-.6,-1) node[above]{$e_1$};
    \draw[very thick,->] (-1,-1) -- (-1,-.6) node[below left]{$e_2$};

    \draw[thick] (-2,-2) rectangle +(2,2);

    \filldraw[fill=gray!30,thick] ( 2,0) rectangle +( 2,2);
    \draw[thick] ( 0,2) rectangle +( 2,2);
    \draw[thick] (-4,0) rectangle +(-2,2);
    
    \node at (0,-1.5) [right] {$X$};
    \node at (4,.5) [right] {$Y_{(2,1)}$};
    \node at (-4,.5) [right] {$Y_{(-2,1)}$};
    \node at (1,2) [below] {$Y_{(1,2)}$};

    \foreach \x/\xt in {-.65/1,.0/2,.65/3} \foreach \y/\yt in
    {-.65/1,.0/2,.65/3} { \fill (\x-5,\y+1) circle (.06)
      node[red,below,font=\tiny] {(\xt,\yt)}; \fill (\x-1,\y-1) circle (.06)
      node[red,below,font=\tiny] {(\xt,\yt)}; \fill (\x+1,\y+3) circle (.06)
      node[red,below,font=\tiny] {(\xt,\yt)}; \fill (\x+3,\y+1) circle (.06)
      node[red,below,font=\tiny] {(\xt,\yt)}; }

    \draw[red,thick,->] (2.2,-.2) to node[below]{$\beta_1$} (3.8,-.2);
    \draw[red,thick,->] (1.8,.2) to node[below left]{$\beta_2$} (1.8,1.8);

    \draw[red,thick,->] (-1.8,-2.2) to node[below right]{$\alpha_1$} (-.2,-2.2);
    \draw[red,thick,->] (-2.2,-1.8) to node[above left]{$\alpha_2$} (-2.2,-.2);

    \fill[white] (-4.2,2.4) rectangle (-1.8,3.6);
    \node at (-3,3) [red,above] {$\alpha = (\alpha_1,\alpha_2)$};
    \node at (-3,3) [red,below] {$\beta = (\beta_1,\beta_2)$};

  \end{tikzpicture}
  \caption{Axial and diagonal symmetries of interactions. The interpolation
    points $\bar x_\alpha$ and $\bar y_\beta$ are indexed by the multi-indices
    $\alpha$ and $\beta$, respectively (interpolation order $\ell=3$). The
    only transfer vector that satisfies $t_1 \ge t_2 \ge 0$ is $t=(2,1)$. In
    that case, we claim that $\Mat{K}_{(2,1)}$ is the only M2L operator we
    need to compute. The direction of the arrows indicates the growth of the
    respective multi-index component.}
  \label{fig:example_axial_diag_sym}
\end{figure}

\paragraph{Model example}
The target cell $X$ in Fig.~\ref{fig:example_axial_diag_sym} has three
interactions $Y_t$ with the transfer vectors $t \in \{(2,1), (1,2),
(-2,1)\}$. We choose the reference domain to be given by $t_1 \ge t_2 \ge
0$. The goal is to express the M2L operators of all interactions via M2L
operators of interactions that lie in the reference domain only. In our
example this is the interaction with the transfer vector $(2,1)$. The two
transfer vectors $(1,2)$ and $(-2,1)$ can be shown to be reflections along the
lines given by $t_1=0$ and $t_1=t_2$, respectively. We easily find that any $t
\in \mathbb{Z}^2$ can be expressed as a reflection (or a combination of
reflections) of transfer vectors that satisfy $t_1 \ge t_2 \ge 0$.

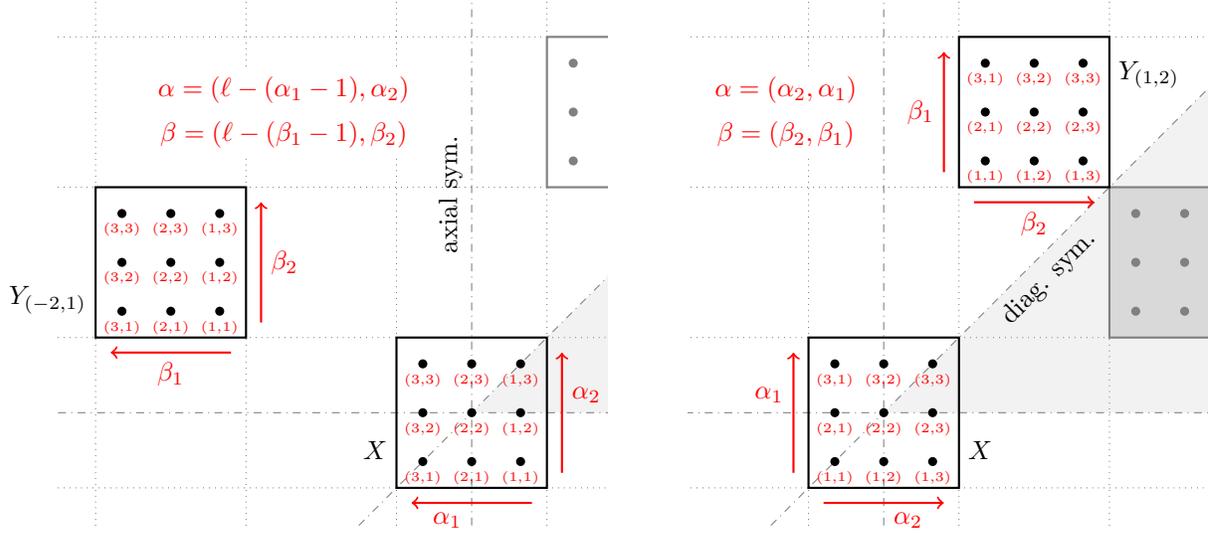
\begin{figure}[h!]
  \centering \subfloat[Invert the component $\alpha_1$ and $\beta_1$ due to
  the axial symmetry $t_1=0$.]{
    \label{fig:example_axial_sym}
    \begin{tikzpicture}
      \clip (-7.5,-2.7) rectangle (.8,4.5);

      \draw[very thin,dotted,step=2] (-6.5,-2.5) grid (4.5,4.5);

      \fill[gray, fill opacity=.1] (-1,-1) -- (4.5,-1) -- (4.5,4.5);

      \draw[gray,dashdotted] (-6.5,-1) -- (4.5,-1);
      \draw[gray,dashdotted] (-1,-2.5) to node[sloped,above right,black]
      {axial sym.} (-1,4.5);
      \draw[gray,dashdotted] (-2.5,-2.5) -- (4.5,4.5);

      \draw[     thick] (-2,-2) rectangle +( 2,2);
      \draw[     thick] (-4, 0) rectangle +(-2,2);
      \draw[gray,thick] ( 0, 2) rectangle +( 2,2);

      \foreach \x/\xt in {-.65/3,.0/2,.65/1}
      \foreach \y/\yt in {-.65/1,.0/2,.65/3}
      {
        \fill (\x-5,\y+1) circle (.06) node[red,below,font=\tiny] {(\xt,\yt)}; 
        \fill (\x-1,\y-1) circle (.06) node[red,below,font=\tiny] {(\xt,\yt)}; 
        \fill[gray] (\x+1,\y+3) circle (.06); 
      }
      
      \node at (-2,-1.5) [left] {$X$};
      \node at (-6,.5) [left] {$Y_{(-2,1)}$};

      \draw[red,thick,<-] (-5.8,-.2) to node[below]{$\beta_1$} (-4.2,-.2);
      \draw[red,thick,->] (-3.8,.2) to node[right]{$\beta_2$} (-3.8,1.8);
      
      \draw[red,thick,<-] (-1.8,-2.2) to node[below left]{$\alpha_1$}
      (-.2,-2.2);
      \draw[red,thick,->] (.2,-1.8) to node[above right]{$\alpha_2$} (.2,-.2);

      \fill[white] (-4.2,2.4) rectangle (-1.8,3.6);
      \node at (-3.5,3) [red,above] {$\alpha = (\ell-(\alpha_1-1),\alpha_2)$};
      \node at (-3.5,3) [red,below] {$\beta  = (\ell-(\beta_1-1), \beta_2)$};

    \end{tikzpicture}
  }
  \hfill
  \subfloat[Swap the components $\alpha_1 \leftrightarrow \alpha_2$ and
  $\beta_1 \leftrightarrow \beta_2$ due to the diagonal symmetry $t_1=t_2$.]{
    \label{fig:example_diag_sym}
    \begin{tikzpicture}
      \clip (-3.6,-2.7) rectangle (3.3,4.5);

      \draw[very thin,dotted,step=2] (-6.5,-2.5) grid (4.5,4.5);

      \fill[gray, fill opacity=.1] (-1,-1) -- (4.5,-1) -- (4.5,4.5);

      \draw[gray,dashdotted] (-6.5,-1) -- (4.5,-1);
      \draw[gray,dashdotted] (-1,-2.5) -- (-1,4.5);
      \draw[gray,dashdotted] (-2.5,-2.5) to node[sloped,below,black]
      {diag. sym.} (4.5,4.5);

      \draw[thick] (-2,-2) rectangle +(2,2);

      \filldraw[gray,fill=gray!30,thick] ( 2,0) rectangle +( 2,2);
      \draw[thick] ( 0,2) rectangle +( 2,2);
      \draw[thick] (-4,0) rectangle +(-2,2);

      \foreach \x/\xt in {-.65/1,.0/2,.65/3}
      \foreach \y/\yt in {-.65/1,.0/2,.65/3}
      {
        \fill (\x-1,\y-1) circle (.06) node[red,below,font=\tiny] {(\yt,\xt)}; 
        \fill (\x+1,\y+3) circle (.06) node[red,below,font=\tiny] {(\yt,\xt)}; 
        \fill[gray] (\x+3,\y+1) circle (.06);
      }
      
      \node at ( 0,-1.5) [right]      {$X$};
      \node at (2,3.5)   [right] {$Y_{(1,2)}$};

      \draw[red,thick,->] (-.2,2.2) to node[left]{$\beta_1$} (-.2,3.8);
      \draw[red,thick,->] (.2,1.8) to node[below]{$\beta_2$} (1.8,1.8);
      
      \draw[red,thick,->] (-1.8,-2.2) to node[below right]{$\alpha_2$}
      (-.2,-2.2);
      \draw[red,thick,->] (-2.2,-1.8) to node[above left]{$\alpha_1$}
      (-2.2,-.2);

      \fill[fill=white] (-3.4,2.4) rectangle (-1.2,3.6);
      \node at (-2.3,3) [red,above] {$\alpha = (\alpha_2,\alpha_1)$};
      \node at (-2.3,3) [red,below] {$\beta = (\beta_2,\beta_1)$};

    \end{tikzpicture}
  }
  \caption{The direction of the arrows indicates the growth of the respective
    multi-index component such that the M2L operators $\Mat{K}_{(-2,1)}$ and
    $\Mat{K}_{(1,2)}$ become the same as $\Mat{K}_{(2,1)}$. In other words,
    the mapping from the arrows in Fig.~\ref{fig:example_axial_diag_sym} to
    the arrows here is analog to the mapping of the multi-indices and
    corresponds to the permutations of $\Mat{K}_{(2,1)}$ such that it
    coincides with $\Mat{K}_{(-2,1)}$ and $\Mat{K}_{(1,2)}$.}
  \label{fig:example_permutations}
\end{figure}

We claim that any reflection of a transfer vector corresponds to a permutation
of the respective M2L operator. Recall that the evaluation of $K(\bar
x_\alpha, \bar y_\beta)$ gives the entry from row $m(\alpha)$ and column
$n(\beta)$ of the M2L operator. $\Mat{K}_{(2,1)}$ is the only M2L operator
whose transfer vector satisfies $t_1 \ge t_2 \ge 0$. The multi-indices are
constructed as presented in Fig.~\ref{fig:example_axial_diag_sym}. As can be
checked, the entry $(\Mat{K}_{(2,1)})_{mn}$ is not the same as
$(\Mat{K}_{(1,2)})_{mn}$ or $(\Mat{K}_{(-2,1)})_{mn}$. However, if we use the
permuted multi-indices from Fig.~\ref{fig:example_axial_sym} for
$\Mat{K}_{(-2,1)}$ or those from Fig.~\ref{fig:example_diag_sym} for
$\Mat{K}_{(1,2)}$ they are the same. The logic behind this can be summarized as
follows.
\begin{itemize}
\item If an \emph{axial symmetry} is given by $t_1=0$ as shown in
  Fig.~\ref{fig:example_axial_sym}, we invert the corresponding component of
  the multi-index as
  \begin{equation}
    \label{eq:axial_multi_index}
    \alpha \leftarrow (\ell-(\alpha_1-1), \alpha_2)
    \quad \text{and} \quad
    \beta \leftarrow (\ell-(\beta_1-1), \beta_2).
  \end{equation}
\item If the \emph{diagonal symmetry} is given by $t_1=t_2$ as shown in
  Fig.~\ref{fig:example_diag_sym}, we swap the corresponding components as
  \begin{equation}
    \label{eq:diag_multi_index}
    \alpha \leftarrow (\alpha_2, \alpha_1)
    \quad \text{and} \quad
    \beta \leftarrow (\beta_2, \beta_1).
  \end{equation}
\end{itemize}
Sometimes it is necessary to combine axial and diagonal permutations. Take as
example the transfer vector~$(-1,2)$: we need to flip it along $t_1=0$ and
then along $t_1=t_2$ to get $(2,1)$. Note that reflections are non
commutative, i.e., the order of their application matters. This is also true
for permutations of the M2L operators.

\paragraph{Generalization}
Let us extend the above concept to the three dimensional case. We start by
introducing three axial and three diagonal symmetries in $\mathbb{Z}^3$.
\begin{itemize}
\item \emph{Axial symmetry planes} are given by $t_1=0$, $t_2=0$ and $t_3=0$
  (see Fig.~\ref{fig:interaction_box}). Each of the three planes divides
  $\mathbb{Z}^3$ in two parts, i.e., the negative part $t_i<0$ and the
  positive part $t_i\ge 0$. By combining all three planes $\mathbb{Z}^3$ is
  divided into octants. In the following we use $\mathbb{Z}_+^3$, i.e., the
  octant with $t_1,t_2,t_3\ge 0$ as reference octant.
  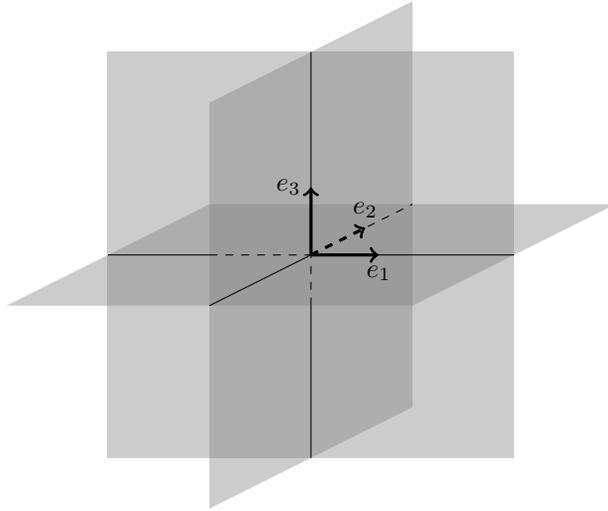
\begin{figure}[htbp]
    \centering
    \begin{tikzpicture}[scale=.9]
      \fill[fill=gray, fill opacity=.4]
      (0,-3) -- (3,-1.5) -- (3,4.5) -- (0,3);

      \fill[fill=gray, fill opacity=.4]
      (-3,0) -- (3,0) -- (6,1.5) -- (0,1.5);

      \fill[fill=gray, fill opacity=.4]
      (-1.5,-2.25) -- (4.5,-2.25) -- (4.5,3.75) -- (-1.5,3.75);

      \draw[thin] (-1.5,.75) -- (0,.75);
      \draw[thin,dashed] (0,.75) -- (1.5,.75);
      \draw[thin] (1.5,.75) -- (4.5,.75);
      \draw[thin] (1.5,-2.25) -- (1.5,0);
      \draw[thin,dashed] (1.5,0) -- (1.5,.75);
      \draw[thin] (1.5,.75) -- (1.5,3.75);
      \draw[thin] (0,0) -- (1.5,.75);
      \draw[thin,dashed] (1.5,.75) -- (3,1.5);

      \draw[very thick,->] (1.5,.75) -- (2.5, .75) node[below]{$e_1$};
      \draw[dashed, very thick,->] (1.5,.75) -- (2.3,1.15) node[above]{$e_2$};
      \draw[very thick,->] (1.5,.75) -- (1.5,1.75) node[left]{$e_3$};
    \end{tikzpicture}
    \caption{Three axial symmetry planes split $\mathbb{Z}^3$ in octants. The
      reference octant is given by $t_1,t_2,t_3\ge 0$.}
    \label{fig:interaction_box}
  \end{figure}
\item \emph{Diagonal symmetry planes} are given by $t_1=t_2$, $t_1=t_3$ and
  $t_2=t_3$ (see Fig.~\ref{fig:symmetries}). In $\mathbb{Z}^3$ there are six
  diagonal symmetries; however, we restrict ourselves to the symmetries
  affecting the reference octant $\mathbb{Z}_+^3$.
  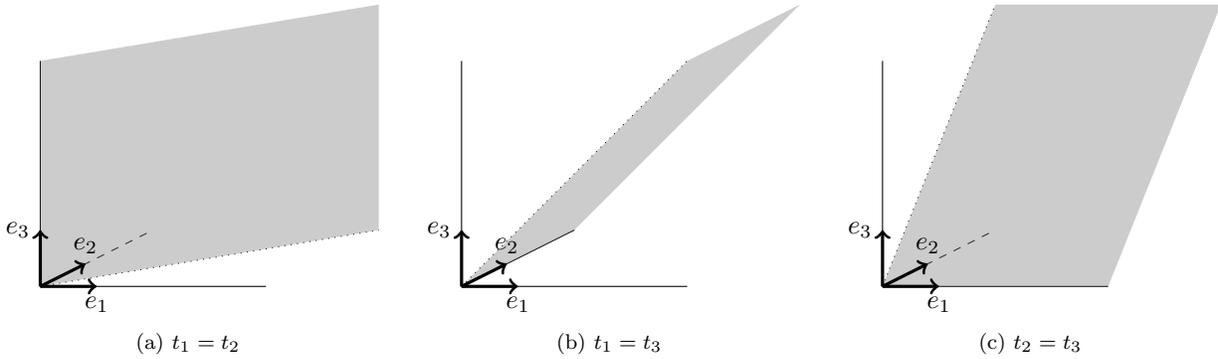
\begin{figure}[htbp]
    \centering
    \subfloat[$t_1 = t_2$]{
      \label{fig:ji}
      \begin{tikzpicture}[scale=.75]
        \draw[thin] (0,4) -- (0,0) -- (4,0);
        \draw[thin,dashed] (0,0) -- (2,1);
        \draw[dotted] (0,0) -- (6,1);

        \fill[fill=gray, fill opacity=.4]
        (0,0) -- (6,1) -- (6,5) -- (0,4);

        \draw[very thick,->] (0,0) -- (1,0)   node[below]{$e_1$};
        \draw[very thick,->] (0,0) -- (.8,.4) node[above]{$e_2$};
        \draw[very thick,->] (0,0) -- (0,1)   node[left]{$e_3$};
      \end{tikzpicture}
    }
    \hfill
    \subfloat[$t_1 = t_3$]{
      \label{fig:ki}
      \begin{tikzpicture}[scale=.75]
        \draw[thin] (0,4) -- (0,0) -- (4,0);
        \draw[thin] (0,0) -- (2,1);
        \draw[dotted] (0,0) -- (4,4);

        \fill[thick, fill=gray, fill opacity=.4]
        (0,0) -- (4,4) -- (6,5) -- (2,1) -- cycle;

        \draw[very thick,->] (0,0) -- (1,0)   node[below]{$e_1$};
        \draw[very thick,->] (0,0) -- (.8,.4) node[above]{$e_2$};
        \draw[very thick,->] (0,0) -- (0,1)   node[left]{$e_3$};
      \end{tikzpicture}
    }
    \hfill
    \subfloat[$t_2 = t_3$]{
      \label{fig:kj}
      \begin{tikzpicture}[scale=.75]
        \draw[thin] (0,4) -- (0,0) -- (4,0);
        \draw[thin,dashed] (0,0) -- (2,1);
        \draw[dotted] (0,0) -- (2,5);

        \fill[fill=gray, fill opacity=.4]
        (0,0) -- (4,0) -- (6,5) -- (2,5) -- cycle;

        \draw[very thick,->] (0,0) -- (1,0)   node[below]{$e_1$};
        \draw[very thick,->] (0,0) -- (.8,.4) node[above]{$e_2$};
        \draw[very thick,->] (0,0) -- (0,1)   node[left]{$e_3$};
      \end{tikzpicture}
    }
    \caption{Three diagonal symmetry planes in the reference octant.}
    \label{fig:symmetries}
  \end{figure}
  \begin{figure}[htbp]
    \centering
    \begin{tikzpicture}[scale=.75]
      \draw[thin] (0,4) -- (0,0) -- (4,0);
      \draw[thin,dashed] (0,0) -- (2,1);
      \fill[gray, opacity=.4] (0,0) -- (6,1) -- (6,5) -- cycle;
      \fill[gray, opacity=.4] (0,0) -- (4,0) -- (6,1) -- cycle;
      \fill[gray, opacity=.4] (0,0) -- (4,0) -- (6,5) -- cycle;
      \draw[dotted] (6,1) -- (0,0) -- (6,5);
      \draw[very thick,->] (0,0) -- (1,0)   node[below]{$e_1$};
      \draw[very thick,->] (0,0) -- (.8,.4) node[above]{$e_2$};
      \draw[very thick,->] (0,0) -- (0,1)   node[left]{$e_3$};
    \end{tikzpicture}
    \caption{The final cone $\mathbb{Z}_{\text{sym}}^3$ ($t_1 \ge t_2 \ge t_3
      \ge 0$) is obtained by combining axial and diagonal symmetries.}
    \label{fig:cone_sym}
  \end{figure}
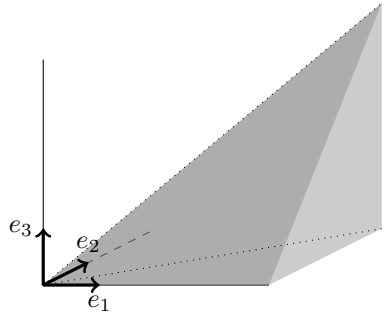
\end{itemize}
By combining the three diagonal symmetries and the three axial symmetries we
obtain the cone shown in Fig.~\ref{fig:cone_sym}. We refer to it as
$\mathbb{Z}_{\text{sym}}^3$; it is given by
\begin{equation}
  \label{eq:z_sym}
  \mathbb{Z}_{\text{sym}}^3 = \left\{\mathbb{Z}_{\text{sym}}^3 \subset \mathbb{Z}^3
  : t_1\ge t_2 \ge t_3 \ge 0 \text{ with } t\in\mathbb{Z}^3\right\}.
\end{equation}
By its means we can identify the subset of transfer vectors $T_{\text{sym}}
\subset T \subset \mathbb{Z}^3$ as
\begin{equation}
  \label{eq:reduced_interaction_list}
  T_{\text{sym}} = T \cap \mathbb{Z}_{\text{sym}}^3
\end{equation}
such that all others $T\backslash T_{\text{sym}}$ can be expressed as
reflections of this fundamental set. Next, we claim that these symmetries are
also useful for M2L operators.

\paragraph{Permutation matrices}
Any reflection of a transfer vector along a symmetry plane determines the
permutation of its associated M2L operator as
\begin{equation}
  \label{eq:permutations}
  \Mat{K}_t = \Mat{P}_t \Mat{K}_{p(t)} \Mat{P}_t^\top.
\end{equation}
The permutation matrix $\Mat{P}_t$ depends on the transfer vector $t\in T$. We
also need the surjective mapping $p : T \rightarrow T_{\text{sym}}$; it
associates every transfer vector in $T$ to exactly one in
$T_{\text{sym}}$. The left application of $\Mat{P}_t$, essentially,
corresponds to the permutation of $\alpha$ and its right application to the
permutation of $\beta$, affecting rows (respectively columns) of the original
matrix $\Mat{K}_{p(t)}$. Note, the permutation matrices $\Mat{P}_t$ depend
only on the transfer vector $t$. How do we construct them? For some $t$ we
introduce axial and diagonal permutations $\pi_t^{\text{A}}$ and
$\pi_t^{\text{D}}$ that read as follows.
\begin{itemize}
\item \emph{Axial symmetries}: multi-index permutations are computed as
  \begin{equation}
    \label{eq:axial_permutation}
    \pi_t^{\text{A}}(\alpha_1, \alpha_2, \alpha_3) = \left(\bar \alpha_1, \bar
      \alpha_2, \bar \alpha_3\right) 
    \quad \text{with} \quad
    \bar \alpha_i =
    \begin{cases}
      \alpha_i & \text{if } t_i \ge 0, \\
      \ell - (\alpha_i-1) & \text{else}.
    \end{cases}
  \end{equation}
  There exist $8$ different possibilities that correspond to the octants
  presented in Fig.~\ref{fig:interaction_box}. Note, $\pi_t^{A}(\alpha) =
  \alpha$ is only true for transfer vectors with $t_1,t_2,t_3 \ge 0$.
\item \emph{Diagonal symmetries}: multi-index permutations are computed as
  \begin{equation}
    \label{eq:diagonal_permutation}
    \pi_t^{\text{D}}(\alpha_1, \alpha_2, \alpha_3) = \left(\alpha_i, \alpha_j,
      \alpha_k\right) \quad \text{such that} \quad |t_i| \ge |t_j| \ge |t_k|.
  \end{equation}
  There exist $6$ different possibilities that correspond to the $6$ different
  cones if we consider Fig.~\ref{fig:cone_sym}. Note again,
  $\pi_t^{\text{D}}(\alpha) = \alpha$ is only true for transfer vectors with
  $t_1 \ge t_2 \ge t_3 \ge 0$.
\end{itemize}
Given these multi-index permutations and the mapping functions $m(\alpha)$ and
$n(\beta)$ we can define a permutation matrix $\Mat{P}_t$ of size $\ell^3
\times \ell^3$. Its entries are $0$ except in column $j$ the entry $i =
m(\pi_t(m^{-1}(j)))$ is $1$. Let us go through the computation of this index:
first, we compute the multi-index $\alpha = m^{-1}(j)$, then, we permute the
multi-index $\bar \alpha = \pi_t(\alpha)$ and last, we compute the row-index
$i = m(\bar \alpha)$. Permutation matrices may be written as
\begin{equation}
  \label{eq:permutation_matrix}
  \Mat{P}_t = \left(e_{m\left(\pi_t(m^{-1}(0))\right)},
    e_{m\left(\pi_t(m^{-1}(1))\right)}, \dots,
    e_{m\left(\pi_t(m^{-1}((\ell-1)^3))\right)}\right), 
\end{equation}
where $e_j$ denotes a column unit vector of length $\ell^3$ with $1$ in the
$j$th position and $0$ elsewhere. Permutation matrices are orthogonal
$\Mat{P}_t \Mat{P}_t^\top = \Mat{I}$, hence, the inverse exists and can be
written as $\Mat{P}_t^{-1} = \Mat{P}_t^\top$. Note that the combination of
permutations is non commutative. Given these permutations $\pi_t^{\text{A}}$
and $\pi_t^{\text{D}}$ we setup $\Mat{P}_t^{\text{A}}$ and
$\Mat{P}_t^{\text{D}}$ and construct the permutation matrix as
\begin{equation}
  \label{eq:final_permutation_matrix}
  \Mat{P}_t = \Mat{P}_t^{\text{D}} \, \Mat{P}_t^{\text{A}}.
\end{equation}
The permutation for the multi-index $\beta$ is the same.

\subsubsection{IA with symmetries (IAsym)}
\label{sec:blocked_convolution}

By exploiting the above introduced symmetries we end up with an optimization
we refer to as the IAsym variant. We individually approximate and store only
M2L operators with $t \in T_{\text{sym}}$ and express all others via
permutations as shown Eqn.~\eqref{eq:permutations}. The IAsym variant for an
arbitrary transfer vector $t\in T$ consists of the following three steps.
\begin{enumerate}
\item Permute multipole expansions
  \begin{equation}
    \label{eq:perm_mult}
    \Mat{w}_t = \Mat{P}_t^\top \Mat{w}
  \end{equation}
\item Compute permuted local expansions
  \begin{equation}
    \label{eq:perm_mult_perm_loc}
    \Mat{f}_t = \Mat{K}_{p(t)} \Mat{w}_t
  \end{equation}
\item Un-permute local expansions
  \begin{equation}
    \label{eq:perm_loc}
    \Mat{f} = \Mat{P}_t \Mat{f}_t
  \end{equation}
\end{enumerate}
Note that the permutation matrix is not applied to the actual M2L operator
(remains unchanged as can be seen in step $2$). Its application is implemented
as a reordering of vector entries (step $1$ and $3$). Depending on whether the
M2L operator exist in its full-rank or in its low-rank representation (see
Eqn.~\eqref{eq:ialowrank}) the application corresponds to one or two
matrix-vector products. In the following we introduce a blocking scheme that
leads to a faster execution on a computer.

\subsubsection{IAsym with blocking (IAblk)}
\label{sec:iasym}

We know from Sec.~\ref{sec:symmetries} that based on the consideration of
symmetries and permutations, many interactions share the same M2L
operators. This paves the road for blocking schemes. Essentially, the idea is
to substitute many matrix-vector products by a few matrix-matrix
products. Blocking schemes do not change the overall complexity of the
algorithm, but they allow for the use of highly optimized matrix-matrix
product implementations. Such achieve much higher peak performances than
optimized matrix-vector product implementations due to better cache reuse and
less memory traffic \citep{Dongarra:1990, intel_mkl}.

\begin{algorithm}
  \caption{Blocking scheme with $|T_{\text{sym}}|$ matrix-matrix products}
  \label{alg:blockedm2l}
  \begin{algorithmic}[1]
    \Function{BlockedM2L}{target cell $X$ and all far-field interactions $I_Y$}
    \State{allocate $\Mat{F}_p$ and $\Mat{W}_p$ for
      $p=1,\dots,|T|_{\text{sym}}$} \label{permmatWF}
    \State{retrieve $\Mat{f}$ from $X$}
    \State{set all $c_p = 0$} \label{expcounter}  
    \For{all source cells $Y$ in $I_Y$}
    \State{retrieve $\Mat{w}$ from $Y$ and compute $t$ from cell-pair $(X,Y)$}
    \State{column $c_{p(t)}$ of $\Mat{W}_{p(t)}$ gets $\Mat{P}_t^\top
      \Mat{w}$ \Comment{Permute multipole expansions}}
    \State{increment $c_{p(t)}$}
    \EndFor
    \For{all $\{\Mat{K}_p\}$}
    \State{$\Mat{F}_p \gets \Mat{K}_p \Mat{W}_p$ \label{matrix_matrix_mult}
      \Comment{Compute permuted local expansions}}
    \EndFor
    \State{set all $c_p = 0$}
    \For{all source cells $Y$ in $I_Y$}
    \State{compute $t$ from cell-pair $(X,Y)$}
    \State{retrieve $\Mat{f}_t$ from column $c_{p(t)}$ of
      $\Mat{F}_{p(t)}$} 
    \State{increment $c_{p(t)}$}
    \State{$\Mat{f} \gets \Mat{f} + \Mat{P}_t \Mat{f}_t$
      \Comment{Permute permuted local expansions}} 
    \EndFor
    \EndFunction
  \end{algorithmic}
\end{algorithm}
In our concrete case, we use the blocking scheme to block multipole and local
expansions. Instead of permuting them and applying the M2L operators
individually (matrix-vector products), we assemble those that share the same
M2L operator as distinct matrices to whom we apply the M2L operators then
(matrix-matrix products). Algorithm~\ref{alg:blockedm2l} explains this in
details. We need the matrices $\Mat{W}_p$ and $\Mat{F}_p$ of size $\ell^2
\times n_p$ for $p=1,\dots,|T_{\text{sym}}|$. Their columns store the permuted
multipole and the resulting (also permuted) local expansions. The values for
$n_p$ indicate how many interactions in $T$ share the same M2L operator of
interactions in $T_{\text{sym}}$, in other words, $n_1 + \dots + n_p + \dots +
n_{|T_{\text{sym}}|} = |T|$ is true. In the case of bbFMM this is a priori
known, since the full interaction list is a priori given (see
Sec.~\ref{sec:translation_invariance}). That is not the case for dFMM and the
values for $n_p$ have to be determined during a precomputation step. We also
need counters $c_p$ to indicate the position of the currently processed
expansions in $\Mat{W}_p$ and $\Mat{F}_p$. As opposed to IAsym, here we split
up the single loop over all interactions into three loops. In the first one,
we assemble the set of matrices $\Mat{W}_p$. At the end $c_p \le n_p$ is true
for all $p$. In the second loop, we perform at most $|T_{\text{sym}}|$
matrix-matrix products. And in the last loop, we increment the local expansion
with the expansions from all $\Mat{F}_p$.

\paragraph{Blocking along multiple target cells}
Algorithm~\ref{alg:blockedm2l} proposes to use the blocking scheme for all
interactions of only \emph{one} target cell. In the worst case no M2L operator
is shared and the algorithm coincides with IAsym. Moreover, the size of the
matrices $\Mat{W}_p$ and $\Mat{F}_p$ might vary to a large extent. That is why
we pursued the blocking idea further to come up with a more efficient
scheme. Instead of using individual $n_p$ we choose it to be a constant $n$
for all $p=1,\dots,|T_{\text{sym}}|$. Then we keep on blocking expansions
using interactions lists of \emph{multiple} (as opposed to one) target
cells. Once $c_p=n$ is true for some $p$, we apply the M2L operator as
$\Mat{F}_p = \Mat{K}_p \Mat{W}_p$ where $\Mat{W}_p, \Mat{F}_p$ are both of
size $\ell^3 \times n$. In our numerical studies we use this blocking scheme
with $n=128$.

\section{Numerical results}
\label{sec:results}

In the previous sections we introduced various optimizations of the M2L
operators for bbFMM and dFMM. As representative kernel functions we use
\begin{equation}
  \label{eq:kernel_functions}
  \text{the Laplace kernel} \quad K(x,y) = \frac{1}{4\pi|x-y|}
  \quad \text{and the Helmholtz kernel} \quad
  K(x,y) = \frac{e^{\imath k|x-y|}}{4\pi|x-y|}.
\end{equation}
In the numerical studies, hereafter, we use the same parameter setting as are
used in the respective publications and for dFMM we use the wavenumber
$k=1$. All computations are performed on a single CPU of a Intel Core
i7--2760QM CPU @ 2.40GHz $\times$ 8 with 8GB shared memory. We used the
compiler gcc~4.6.3 with the flags ``-O2 -ffast-math''. The source files (C++)
of the program we used can be downloaded via
\verb!http://github.com/burgerbua/dfmm!; they are distributed under the ``BSD
2-Clause'' license.

\paragraph{M2L optimizations}
We show and analyze results for the following eight variants:
\begin{enumerate}
\item NA: the full interaction list $T$ is represented by full-rank (not
  approximated) M2L operators
\item NAsym: the reduced interaction list $T_{\text{sym}}$ is represented by
  full-rank (not approximated) M2L operators
\item NAblk: same as NAsym but with additional blocking of multipole and local
  expansions
\item SA: the variant presented in \cite{fong09a} and briefly sketched in
  Sec.~\ref{sec:fast_convolution_former}
\item SArcmp: same as SA but with additional recompression of all $\Mat{C}_t$
\item IA: same as NA but with low-rank (individually approximated) M2L
  operators
\item IAsym: same as NAsym but with low-rank (individually approximated) M2L
  operators
\item IAblk: same as NAblk but with low-rank (individually approximated) M2L
  operators (see Alg.~\ref{alg:blockedm2l})
\end{enumerate}
Moreover, we study two different low-rank approximation schemes for the SA and
IA variants: on one hand we use a truncated singular value decomposition (SVD)
and on the other hand the adaptive cross approximation (ACA) followed by a
truncated SVD \citep{bebendorf:05}.

\begin{figure}[htbp]
  \centering
  \subfloat[sphere]{\label{fig:e0}
    \includegraphics[width=0.3\textwidth]{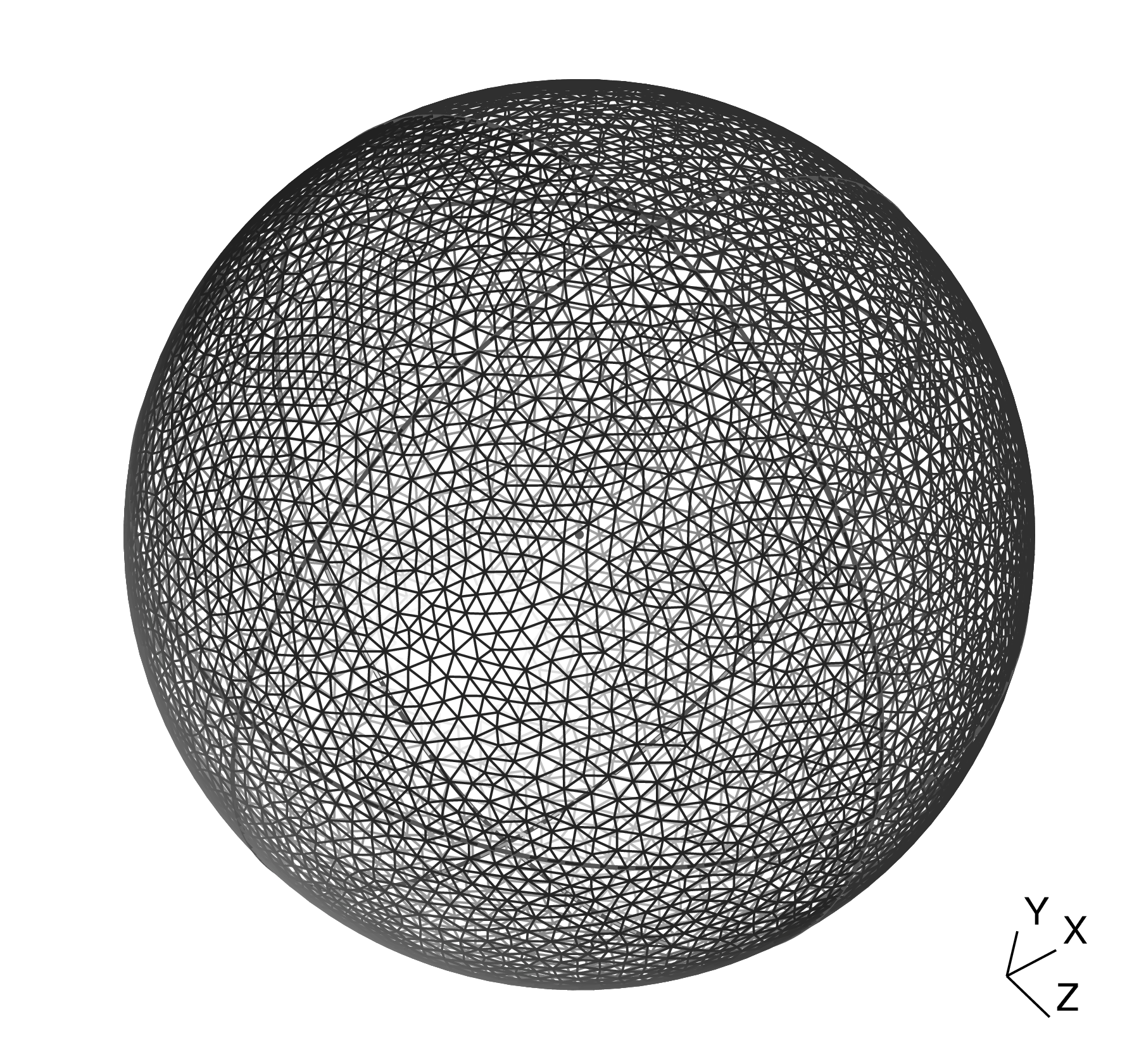}}
  \subfloat[oblate spheroid]{\label{fig:e1}
    \includegraphics[width=0.3\textwidth]{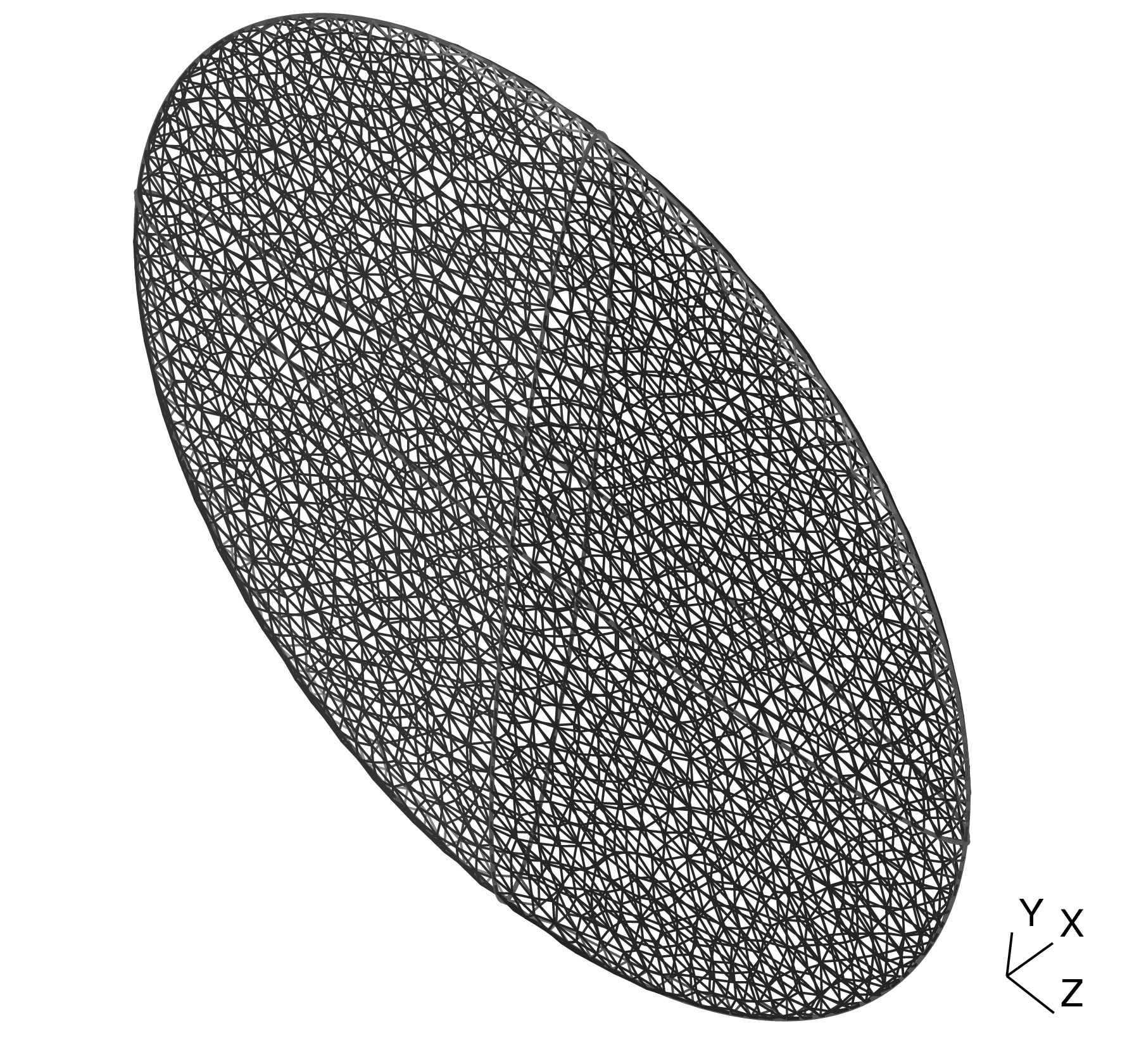}}
  \subfloat[prolate spheroid]{\label{fig:e2}
    \includegraphics[width=0.3\textwidth]{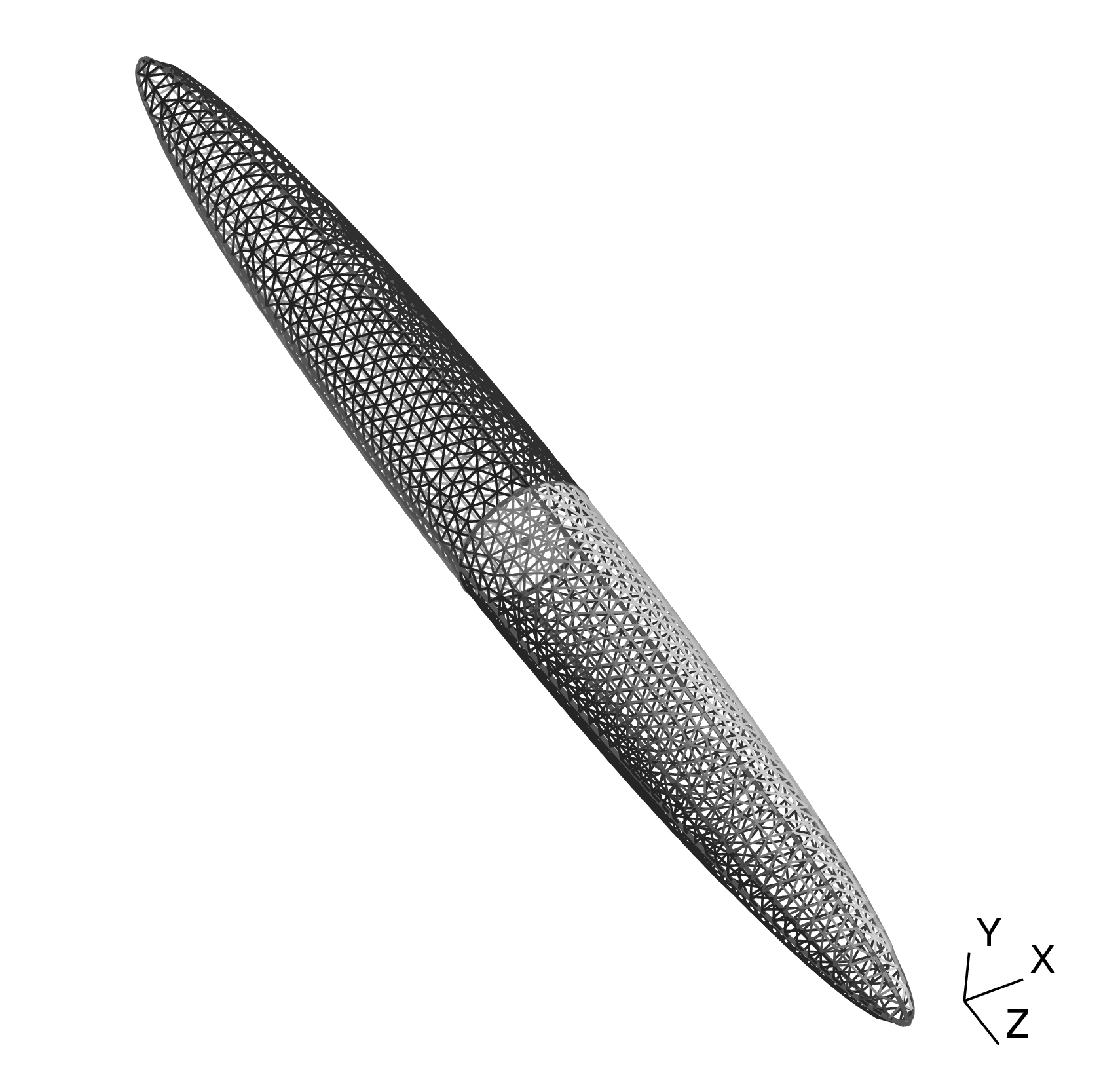}}
  \caption{The example geometries are centered at $(0,0,0)$}
  \label{fig:x_dimensional_objects}
\end{figure}

\paragraph{Example geometries}
We use the three benchmark examples shown in
Fig.~\ref{fig:x_dimensional_objects} and described in the listing below. The
depth of the octree is chosen such that near- and far-field are balanced,
i.e., we want the fastest possible matrix-vector product.
\begin{enumerate}
\item The \emph{sphere} from Fig.~\ref{fig:e0} is contained in the
  bounding-box $64 \times 64 \times 64$; $168\,931$ particles are randomly
  scattered on its surface. The octree has $6$ levels.
\item The \emph{oblate sphere} from Fig.~\ref{fig:e1} is contained in the
  bounding-box $6.4\times64\times64$; $125\,931$ particles are randomly
  scattered on its surface. The octree has $6$ levels.
\item The \emph{prolate sphere} from Fig.~\ref{fig:e2} is contained in the
  bounding-box $6.4\times6.4\times64$; $119\,698$ particles are randomly
  scattered on its surface. The octree has $7$ levels.
\end{enumerate}
In approximative sense, the sphere is a three-dimensional, the oblate sphere a
two-dimensional and the prolate sphere a one-dimensional object in
$\mathbb{R}^3$. We choose these three geometries to study their influence on
the performance on the dFMM. Table~\ref{tab:size_T_all} shows the size of the
full and the reduced interaction lists ($|T|$ and $|T_{\text{sym}}|$) per
level (lf stands for low-frequency, hf for high-frequency regime) for all
three geometries. The size of the interaction lists clearly grows with the
dimensionality of the geometry. We report on the impact of this behavior
later.
\begin{table}[htbp]
  \centering
  \begin{tabular}{c | rrrr | rrrr | rrrr}
    \toprule
    &\multicolumn{4}{c|}{sphere}
    &\multicolumn{4}{c|}{oblate sphere}
    &\multicolumn{4}{c}{prolate sphere} \\
    & 3(hf) & 4(hf) & 5(hf) & 6(hf)
    & 3(hf) & 4(hf) & 5(hf) & 6(hf)
    & 4(hf) & 5(hf) & 6(hf) & 7(lf) \\
    \midrule
    $|T|$            &668&18710&2666&418 &60&2336&1502&400 &214&738&382&424\\
    $|T_{\text{sym}}|$ & 20&  518&  93& 21 & 6& 203&  89& 21 & 35& 61& 21& 22\\
    \bottomrule
  \end{tabular}
  \caption{Size of interactions lists per level for dFMM (hf and lf stands for
    high- and low-frequency regime, respectively)} 
  \label{tab:size_T_all}
\end{table}

\subsection{Accuracy of the method}
\label{sec:accuracy}

\begin{figure}[htbp]
  \centering
  \subfloat[bbFMM for the smooth kernel function]{
    \label{fig:smooth}
    \begin{tikzpicture}
      \begin{loglogaxis}[
        ylabel = Relative error $\varepsilon_{L_2}$,
        xlabel = Target accuracy $\varepsilon_{\text{ACA}}$,
        legend style={at={(0.02,.98)}, anchor=north west},
        ymin = 6e-10,
        ymax = 3e-1
        ]
        \addplot coordinates {
          (1e-1,0.07897)
          (1e-2,0.00534158)
          (1e-3,0.000579394)
          (1e-4,0.000591432)
          (1e-5,0.000111352)
        };
        \addplot coordinates {
          (1e-2,0.0068437)
          (1e-3,0.00520437)
          (1e-4,0.00503505)
          (1e-5,6.07907e-05)
          (1e-6,5.99887e-05)
        };
        \addplot coordinates {
          (1e-3,0.000482389)
          (1e-4,0.000130543)
          (1e-5,8.39802e-06)
          (1e-6,9.91877e-06)
          (1e-7,1.00733e-05)
        };
        \addplot coordinates {
          (1e-4,1.33369e-05)
          (1e-5,4.94247e-07)
          (1e-6,1.27856e-06)
          (1e-7,1.59813e-06)
          (1e-8,1.59031e-06)
        };
        \addplot coordinates {
          (1e-5,1.09018e-06)
          (1e-6,4.36319e-07)
          (1e-7,3.36201e-07)
          (1e-8,2.4486e-07)
          (1e-9,2.45463e-07)
        };
        \addplot coordinates {
          (1e-6,1.23553e-07)
          (1e-7,2.00951e-08)
          (1e-8,3.27886e-08)
          (1e-9,3.7549e-08)
          (1e-10,3.75819e-08)
        };
        \addplot coordinates {
          (1e-7,1.34688e-08)
          (1e-8,3.5816e-09)
          (1e-9,4.7607e-09)
          (1e-10,4.82979e-09)
          (1e-11,4.8339e-09)
        };
        \legend{$\ell=3$, $\ell=4$, $\ell=5$,
          $\ell=6$, $\ell=7$, $\ell=8$, $\ell=9$}
      \end{loglogaxis}
    \end{tikzpicture}}
  \subfloat[dFMM for the oscillatory kernel function]{
    \label{fig:oscillatory}
    \begin{tikzpicture}
      \begin{loglogaxis}[
        xlabel = Target accuracy $\varepsilon_{\text{ACA}}$,
        legend style={at={(0.02,.98)}, anchor=north west},
        ymin = 6e-10,
        ymax = 3e-1
        ]
        \addplot coordinates {
          (1e-1,0.0121434)
          (1e-2,0.0035947)
          (1e-3,0.000828891)
          (1e-4,0.000771916)
          (1e-5,0.000759817)
        };
        \addplot coordinates {
          (1e-2,0.0044148)
          (1e-3,0.000373277)
          (1e-4,0.000125225)
          (1e-5,0.000120835)
          (1e-6,0.000120641)
        };
        \addplot coordinates {
          (1e-3,0.000180832)
          (1e-4,2.62309e-05)
          (1e-5,1.7036e-05)
          (1e-6,1.57456e-05)
          (1e-7,1.58442e-05)
        };
        \addplot coordinates {
          (1e-4,7.02157e-05)
          (1e-5,3.0166e-06)
          (1e-6,1.97545e-06)
          (1e-7,1.86837e-06)
          (1e-8,1.86337e-06)
        };
        \addplot coordinates {
          (1e-5,1.73552e-06)
          (1e-6,5.9981e-07)
          (1e-7,5.11144e-07)
          (1e-8,4.82689e-07)
          (1e-9,4.87089e-07)
        };
        \addplot coordinates {
          (1e-6,1.59941e-07)
          (1e-7,2.21211e-08)
          (1e-8,1.66868e-08)
          (1e-9,1.58526e-08)
          (1e-10,1.58908e-08)
        };
        \addplot coordinates {
          (1e-7,2.3743e-08)
          (1e-8,5.57197e-09)
          (1e-9,6.09465e-09)
          (1e-10,5.98436e-09)
          (1e-11,5.98272e-09)
        };
        \legend{$\ell=3$, $\ell=4$, $\ell=5$,
          $\ell=6$, $\ell=7$, $\ell=8$, $\ell=9$}
      \end{loglogaxis}
    \end{tikzpicture}}
  \caption{Accuracies for the prolate sphere from Fig.~\ref{fig:e2}. The
    target accuracy $\varepsilon_{\text{ACA}}$ refers to the accuracy of the
    approximate M2L operators (see Eqn.~\eqref{eq:ialowrank}). Here, we used
    the ACA followed by a truncated SVD.}
  \label{fig:accuracies}
\end{figure}
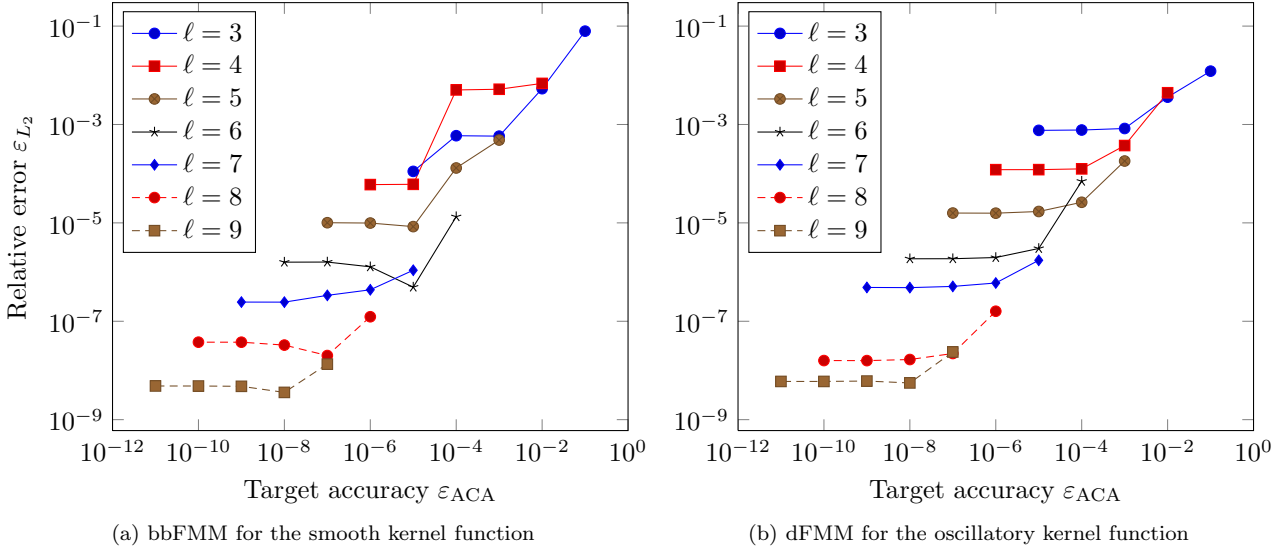

Both, the bbFMM and the dFMM, have two approximations: 1) the interpolation of
the kernel functions determined by interpolation order $\ell$, and 2) the
subsequent low-rank approximation of the M2L operators determined by the
target accuracy $\varepsilon$. The final relative error is a result of both
approximations. We compute it as
\begin{equation}
  \label{eq:rel_err}
  \varepsilon_{L_2} = \left(\frac{\sum_{i\in M} |f_i-\bar f_i|^2}{\sum_{i\in
        M}|f_i|}\right)^{1/2}
\end{equation}
where $M$ is the number of particles $x$ in an arbitrary reference cluster at
the leaf level; $f$ and $\bar f$ are the exact and approximate results,
respectively. In Fig.~\ref{fig:accuracies} we compare achieved accuracies for
the bbFMM and the dFMM with the IAblk variant (other variants produce
identical results). Both plots show the behavior of the relative error
$\varepsilon_{L_2}$ depending on the interpolation order $\ell$ and the target
accuracy $\varepsilon$. Evident is the matching result between the left and
right figure. All curves show an initial plateau and then, after a sharp knee
a slope of roughly $1$. The knee occurs approximately at
$(\ell,\varepsilon)=(Acc, 10^{-Acc})$. In the rest of the paper we use this
convention to describe the accuracy $Acc$ of bbFMM and dFMM.
The low-rank approximations for the computations whose accuracies are shown in
Fig.~\ref{fig:accuracies} were conducted with the ACA followed by a truncated
SVD. By just using a truncated SVD we obtain identical results.

\subsection{Reducing the cost with the SArcmp variant}
\label{sec:recompression_ucb}

The cost of applying an approximate M2L operator mainly depends on its rank
$k$. In Tab.~\ref{tab:lowranks} we compare the average rank of M2L operators
for the SA and the IA variants at all levels that have expansions. Let us
explain how we computed the average rank for SA. Recall, when we use that
variant, all M2L operators from one interaction list posses the same rank; the
bbFMM and the dFMM in the low-frequency regime have one and the dFMM in the
high-frequency regime has potentially multiple (directional) interaction lists
per level. Thus, the average rank per level is the average of all ranks used
in all interactions at that level.

The application of one M2L operator from SA and IA requires $\op{O}(k^2)$,
respectively $\op{O}(2k\ell^3)$ operations. Note, the ranks $k$ of the M2L
operators are different; for a given accuracy $Acc$ they are normally
significantly lower for IA than for SA. This can be seen in
Tab.~\ref{tab:lowranks}.
\begin{table}[htbp]
  \centering
    \begin{tabular}{c | rrrr | rrrr}
    \toprule
    & \multicolumn{4}{c}{SA}
    & \multicolumn{4}{|c}{IA}\\
    $Acc$
    &$4$(hf) & $5$(hf) & $6$(hf) &$7$(lf)
    &$4$(hf) & $5$(hf) & $6$(hf) &$7$(lf)\\
    \midrule
    $3$& 9.8&12.3& 12.8& 19& 5.4& 5.7& 5.7& 5.0\\
    $5$&21.7&30.8& 39.2& 71&11.3&12.3&13.5&12.6\\
    $7$&38.2&58.6& 80.0&163&18.9&21.5&24.7&24.1\\
    $9$&57.8&96.5&138.7&296&28.6&33.2&39.7&40.1\\
    \bottomrule
  \end{tabular}
  \caption{Comparison of average ranks $k$ for the SA and IA variants of
    the dFMM (prolate sphere)}
  \label{tab:lowranks}
\end{table}
The large ranks at level $7$ (first level in the low-frequency regime) of SA
are noteworthy. There, the lower bound for the separation criterion is given
by the usual low-frequency criterion saying that non touching cells are
admissible. Hence, the smallest possible transfer vectors have a length of
$\min_{t\in T} |t| = 2$. The slowly decaying singular values of associated M2L
operators are responsible for the large ranks. On the other hand, the upper
bound for the separation criterion coincides with the lower bound of level $6$
(parent level), which is in the high-frequency regime. Hence, the largest
possible transfer vectors have a length of $\max_{t\in T} |t| \sim 4k$. M2L
operators whose transfer vectors are in that range have much faster decaying
singular values. This fact explains the efficiency of SArcmp (individual
recompression of each M2L operator).

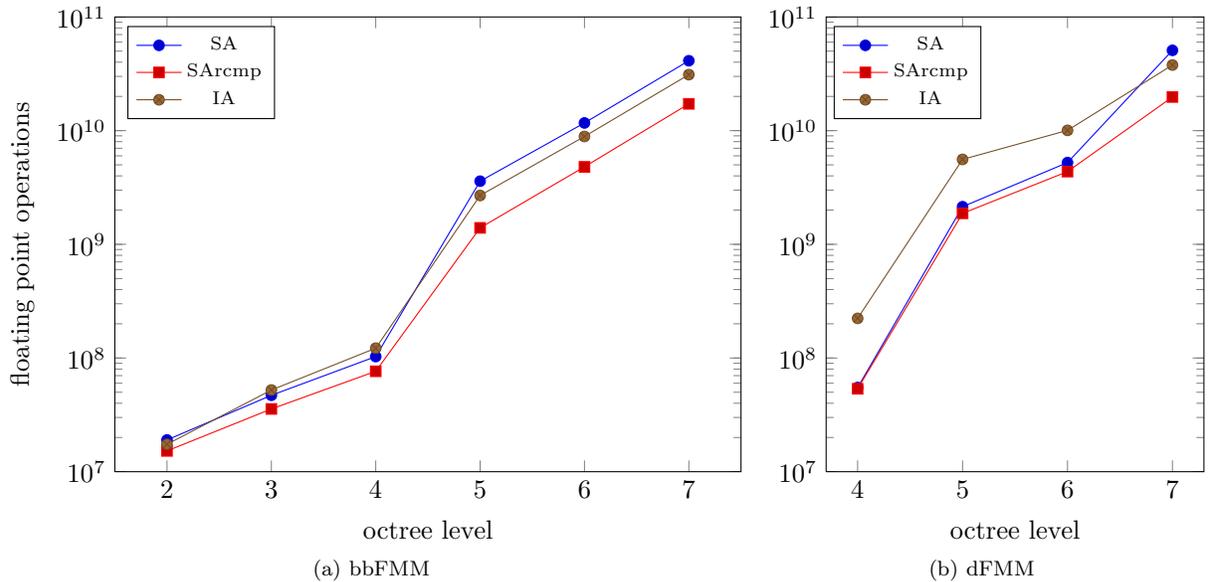
\begin{figure}[htbp]
  \centering
  \subfloat[bbFMM]{\label{fig:flops_per_level_bbfmm}
    \begin{tikzpicture}
      \begin{semilogyaxis}[
        ylabel = floating point operations,
        xlabel = octree level,
        legend style={at={(.02,.98)}, font=\scriptsize, anchor=north west},
        xtick = {2,3,4,5,6,7},
        width = .6\textwidth,
        height = .35\textheight,
        ymin = 1e7,
        ymax = 1e11
        ]
        \addplot coordinates {
          (2,19023984)
          (3,47005920)
          (4,102969792)
          (5,3593159352)
          (6,11685706704)
          (7,41187726531)
        };
        \addplot coordinates {
          (2,15237128)
          (3,35645352)
          (4,76461800)
          (5,1395722080)
          (6,4805612232)
          (7,17193582235)
        };
        \addplot coordinates {
          (2,17431700)
          (3,52295100)
          (4,122021900)
          (5,2689373170)
          (6,8885514660)
          (7,31080557860)
        };
        \legend{SA, SArcmp, IA}  
      \end{semilogyaxis}
    \end{tikzpicture}
  }
  \subfloat[dFMM]{\label{fig:flops_per_level_dfmm}
    \begin{tikzpicture}
      \begin{semilogyaxis}[
        xlabel = octree level,
        legend style={at={(.02,.98)}, font=\scriptsize, anchor=north west},
        xtick = {4,5,6,7},
        width = .4\textwidth,
        height = .35\textheight,
        ymin = 1e7,
        ymax = 1e11
        ]
        \addplot coordinates {
          (4,55016178)
          (5,2142392412)
          (6,5237634064)
          (7,50842136079)
        };
        \addplot coordinates {
          (4,53638474)
          (5,1872228170)
          (6,4362105260)
          (7,19760660025)
        };
        \addplot coordinates {
          (4,223210295)
          (5,5582379495)
          (6,10053526010)
          (7,37764384680)
        };
        \legend{SA, SArcmp, IA}  
      \end{semilogyaxis}
    \end{tikzpicture}
  }
  \caption{Comparison of the M2L cost (floating point operations) growth per
    level for the SA, SArcmp and IA variants of bbFMM and dFMM ($Acc=9$ and
    prolate sphere)}
  \label{fig:flops_per_level}
\end{figure}

In Tab.~\ref{tab:flops_ucbe2} we analyze the cost savings of SArcmp compared
to SA. The left values in each column multiplied by $10^6$ give the overall
number of floating point operations per level for SA. The right values (in
brackets) indicate the respective cost ratio of SArcmp to SA. The
recompression reduces the cost remarkably (see also
Fig.~\ref{fig:flops_per_level}). At the low-frequency level $7$, SArcmp
reduces the cost by more than a factor of $2$. This is almost twice as much as
in high frequency levels. For the impact on timing results we refer to
Sec.~\ref{sec:overalltimings}.
\begin{table}[htbp]
  \centering
    \begin{tabular}{c | rl|rl|rl|rl}
    \toprule
    & \multicolumn{8}{|c}{$\operatorname{cost}(\text{SA}) / 10^6 \quad
      \left(\operatorname{cost}(\text{SArcmp}) /
        \operatorname{cost}(\text{SA})\right)$} 
    \\
    $Acc$
    & \multicolumn{2}{c|}{4(hf)}
    & \multicolumn{2}{c|}{5(hf)}
    & \multicolumn{2}{c|}{6(hf)}
    & \multicolumn{2}{c}{7(lf)} \\
    \midrule
    $3$& 0.8&(0.98)&  20.3&(0.93)&  31.5&(0.93)&  204.3&(0.62)\\
    $5$& 4.9&(0.97)& 161.4&(0.89)& 339.3&(0.86)& 2889.3&(0.47)\\
    $7$&19.3&(0.97)& 687.8&(0.88)&1590.2&(0.83)&15420.0&(0.40)\\
    $9$&55.0&(0.97)&2138.2&(0.87)&5237.6&(0.83)&51505.7&(0.38)\\
    \bottomrule
  \end{tabular}
  \caption{Comparison of cost in terms of floating point operations for SA and
    SArcmp (prolate sphere)} 
  \label{tab:flops_ucbe2}
\end{table}

In Fig.~\ref{fig:flops_per_level} we compare the cost of SA, SArcmp and IA for
bbFMM and dFMM for the prolate sphere and accuracy $Acc=9$. The bbFMM in
Fig.~\ref{fig:flops_per_level_bbfmm} has expansions at the levels $2-7$. The
reason for the jump between level $4$ and $5$ is because the levels $2-4$ have
a maximum of $|T|=16$ and the levels $5-7$ a maximum of $|T|=316$ (common
maximum size for bbFMM) possible M2L operators. The jump in the case of dFMM
in Fig.~\ref{fig:flops_per_level_dfmm} at level $5$ can be explained in the
same way: if we look at Tab.~\ref{tab:size_T_all} we see that $|T|$ is about
twice as large as it is at the other levels. The levels $4,5,6$ of dFMM are in
the high-frequency regime. There, the cost of IA is approximately $12,5,3$
times greater than the cost of SA. However, at level $7$ (low-frequency
regime) of dFMM and at the levels $5-7$ of bbFMM the cost of IA is only about
$2/3$ compared to SA. The reason is the size of the interaction lists
$|T|$. The larger they become the larger the span between the smallest and the
largest rank and that favors individual approximations. The SArcmp is
computationally the least expensive variant. Similar results are obtained for
all other accuracies.

\subsection{Speeding up the precomputation with IA variants}
\label{sec:precomputation}

The bottleneck of SA and SArcmp is the relatively large precomputation time.
This is a minor issue in the case of homogeneous kernel functions. In that
case the approximated M2L operators can be stored on the disk and loaded if
needed for further computations. However, if we deal with non-homogeneous
kernel functions, such as the Helmholtz kernel, IAsym is the way to go. In
Tab.~\ref{tab:prectimings} we compare the precomputation time of SA, IA and
IAsym (we do not report on SArcmp; due to the additional recompression its
precomputation time is higher than the one for SA). For the low-rank
approximation we use a truncated SVD or the ACA followed by a truncated
SVD. In both cases we end up with the same rank. We get remarkable speedups in
the precomputation. Let us look at an extreme example: for the sphere and an
accuracy $Acc=7$ we have $t_{\text{SVD}} = \unit[4740.9]{s}$ for SA versus
$t_{\text{ACA}} = \unit[1.8]{s}$ for IAsym. This corresponds to a speedup
greater than $2600$. For the oblate and prolate sphere we obtain similar results.
\begin{table}[h!]
  \centering
    \begin{tabular}{c|rr|rr|rr}
    \toprule
    & \multicolumn{2}{c}{SA}
    & \multicolumn{2}{|c}{IA}
    & \multicolumn{2}{|c}{IAsym} \\
    $Acc$
    &$\unit[t_{\text{SVD}}]{[s]}$
    &$\unit[t_{\text{ACA}}]{[s]}$
    &$\unit[t_{\text{SVD}}]{[s]}$
    &$\unit[t_{\text{ACA}}]{[s]}$
    &$\unit[t_{\text{SVD}}]{[s]}$
    &$\unit[t_{\text{ACA}}]{[s]}$ \\
    \midrule
    \multicolumn{7}{c}{sphere} \\
    \midrule
    $3$ &   7.0&  4.6&  6.2& 1.9& 0.2&0.1 \\
    $4$ &  75.0& 20.2& 37.0& 4.9& 1.1&0.2 \\
    $5$ & 317.2& 69.1&188.8&12.4& 5.6&0.4 \\
    $6$ &1336.4&197.0&790.5&28.6&22.9&0.9 \\
    $7$ &4740.9&435.9&    -&   -&84.0&1.8 \\
    \midrule
    \multicolumn{7}{c}{oblate sphere} \\
    \midrule
    $3$ &   1.4&  0.9&   1.2& 0.4&  0.1&0.0 \\
    $4$ &  13.4&  3.8&   7.1& 1.1&  0.5&0.1 \\
    $5$ &  59.7& 13.5&  37.0& 2.7&  2.8&0.2 \\
    $6$ & 299.4& 39.4& 150.1& 6.4& 11.1&0.5 \\
    $7$ &1021.2& 97.6& 551.0&13.3& 40.9&0.9 \\
    $8$ &     -&217.7&1751.4&29.0&129.3&2.0 \\
    $9$ &     -&444.6&     -&   -&369.2&3.9 \\
    \midrule
    \multicolumn{7}{c}{prolate sphere} \\
    \midrule
    $3$ &   0.6&   0.7&   0.8 &  0.3 &   0.1 &0.0 \\
    $4$ &   7.5&   4.1&   5.3 &  0.7 &   0.4 &0.1 \\
    $5$ &  64.7&  18.3&  31.0 &  2.1 &   2.5 &0.1 \\
    $6$ & 358.2&  73,3& 199.4 &  5.0 &  15.7 &0.4 \\
    $7$ &1374.8& 204.3& 808.9 & 11.6 &  66.0 &0.8 \\
    $8$ &     -& 549.3&    -  & 25.0 & 322.3 &1.7 \\
    $9$ &     -&1273.2&    -  & 50.2 & 807.4 &3.7 \\
    \bottomrule
  \end{tabular}
  \caption{Precomputation times (SVD versus ACA) for the SA, IA and IAsym
    variants. Missing numbers mean that the available memory of $\unit[8]{GB}$
    has not been sufficient for the respective computation.}  
  \label{tab:prectimings}
\end{table}

\subsection{Comparison of all M2L variants}
\label{sec:overalltimings}

In the previous sections we have revealed the impact of SArcmp and IAsym on the computational cost and the precomputation time of the M2L operators, respectively. In this section we compare all variants and focus on their impact on memory requirement and running time (of applying M2L operators during the actual matrix-vector product). Tables~\ref{tab:timebbfmm_e0}--\ref{tab:timebbfmm_e2} (respectively, \ref{tab:timem2l_e0}--\ref{tab:timem2l_e2}) present running times of the M2L operators for bbFMM (respectively, dFMM) and all three geometries. The upper set of rows reports on results obtained with a BLAS and LAPACK implementation (libblas 1.2.20110419-2, which is sub-optimal in our case). The lower set shows the times obtained with the Intel MKL \citep{intel_mkl} (Intel MKL 10.3.11.339 (intel64)), which proved faster for our purpose.

\begin{table}[htbp]
  \centering
    \begin{tabular}{c|ccc|cc|ccc}
    \toprule
    $Acc$ & NA & NAsym & NAblk & SA & SArcmp & IA & IAsym & IAblk \\
    \midrule
    \multicolumn{9}{c}{libblas 1.2.20110419-2} \\
    \midrule
    3 &  0.7&  0.8&  0.4& 0.5&\textbf{ 0.3}& 0.4& 0.4&\textbf{ 0.3}\\
    4 &  3.6&  3.8&  1.8& 1.4&\textbf{ 0.9}& 1.3& 1.4&\textbf{ 0.7}\\
    5 & 14.2& 12.8&  5.9& 4.6&\textbf{ 2.2}& 3.4& 3.6&\textbf{ 1.8}\\
    6 & 42.8& 38.2& 16.5&11.5&\textbf{ 4.8}& 9.0& 8.7&\textbf{ 4.2}\\
    7 &102.7&101.7& 40.7&25.5&\textbf{ 9.6}&20.2&18.4&\textbf{ 8.6}\\
    8 &229.3&234.2& 89.7&47.8&\textbf{18.3}&40.7&35.9&\textbf{16.2}\\
    9 &    -&484.5&180.0&83.8&\textbf{30.0}&74.0&68.6&\textbf{28.4}\\
    \midrule
    \multicolumn{9}{c}{Intel MKL 10.3.11.339 (intel64)} \\
    \midrule
    3 &  0.3&  0.3& 0.2& 0.2& 0.4& 0.4& 0.4&\textbf{0.2}\\
    4 &  1.8&  1.2& 0.6& 0.6& 0.7& 0.7& 0.8&\textbf{0.4}\\
    5 &  9.4&  6.3& 1.9& 2.0& 1.4& 2.0& 2.0&\textbf{1.0}\\
    6 & 28.2& 14.1& 4.8& 6.9& 2.6& 5.0& 4.0&\textbf{1.7}\\
    7 & 72.6& 57.7&12.0&19.3& 5.8&12.7& 8.2&\textbf{3.5}\\
    8 &127.6&117.0&25.6&34.3&12.0&24.6&16.2&\textbf{6.0}\\
    9 &    -&260.5&50.8&60.4&20.6&44.6&32.9&\textbf{9.9}\\
    \bottomrule
  \end{tabular}
  \caption{M2L timings for the bbFMM (sphere). In this table and below as
    well, bold numbers correspond to the smallest entry in a row. In some
    cases, two columns use bold font when the running times are sufficiently close
    that the difference is not significant.} 
  \label{tab:timebbfmm_e0}
\end{table}
\begin{table}[htbp]
  \centering
    \begin{tabular}{c|ccc|cc|ccc}
    \toprule
    $Acc$ & NA & NAsym & NAblk & SA & SArcmp & IA & IAsym & IAblk \\
    \midrule
    \multicolumn{9}{c}{libblas 1.2.20110419-2} \\
    \midrule
    3 &  0.3&  0.4& 0.2& 0.2&\textbf{ 0.1}& 0.2& 0.2&\textbf{ 0.1}\\
    4 &  1.6&  1.6& 0.8& 0.8&\textbf{ 0.4}& 0.6& 0.6&\textbf{ 0.3}\\
    5 &  6.4&  5.8& 2.6& 1.9&\textbf{ 0.9}& 1.6& 1.6&\textbf{ 0.8}\\
    6 & 19.1& 16.7& 7.4& 5.4&\textbf{ 2.1}& 4.0& 3.9&\textbf{ 1.9}\\
    7 & 46.8& 44.8&18.3&11.2&\textbf{ 4.5}& 9.0& 9.0&\textbf{ 3.8}\\
    8 &102.8&101.9&40.1&21.3&\textbf{ 8.4}&18.2&16.5&\textbf{ 7.2}\\
    9 &    -&212.2&81.1&37.0&\textbf{13.5}&32.4&30.4&\textbf{12.8}\\
    \midrule
    \multicolumn{9}{c}{Intel MKL 10.3.11.339 (intel64)} \\
    \midrule
    3 &  0.1&  0.1& 0.1& 0.1&0.2& 0.2& 0.2&\textbf{0.1}\\
    4 &  0.7&  0.5& 0.3& 0.3&0.3& 0.3& 0.4&\textbf{0.2}\\
    5 &  4.5&  2.0& 0.9& 1.0&0.7& 0.9& 0.9&\textbf{0.4}\\
    6 & 13.1&  6.7& 2.2& 3.2&1.2& 2.6& 1.8&\textbf{0.8}\\
    7 & 33.7& 27.4& 5.4& 8.2&2.7& 5.7& 4.4&\textbf{1.6}\\
    8 & 57.9& 52.7&11.5&14.6&5.2&10.9& 7.2&\textbf{2.7}\\
    9 &117.4&118.2&22.9&27.3&9.1&20.0&14.2&\textbf{4.4}\\
    \bottomrule
  \end{tabular}
  \caption{M2L timings for bbFMM (oblate sphere)}  
  \label{tab:timebbfmm_e1}
\end{table}
\begin{table}[htbp]
  \centering
    \begin{tabular}{c|ccc|cc|ccc}
    \toprule
    $Acc$ & NA & NAsym & NAblk & SA & SArcmp & IA & IAsym & IAblk \\
    \midrule
    \multicolumn{9}{c}{libblas 1.2.20110419-2} \\
    \midrule
    3 &  0.2&  0.2& 0.1& 0.1&\textbf{0.1}& 0.1& 0.1&\textbf{0.1}\\
    4 &  1.0&  1.1& 0.5& 0.4&\textbf{0.3}& 0.4& 0.4&\textbf{0.2}\\
    5 &  4.0&  4.0& 1.9& 1.2&\textbf{0.6}& 1.2& 1.0&\textbf{0.6}\\
    6 & 12.0& 11.8& 4.7& 3.4&\textbf{1.3}& 2.6& 2.4&\textbf{1.4}\\
    7 & 29.3& 31.4&11.7& 6.9&\textbf{2.8}& 7.2& 5.2&\textbf{2.9}\\
    8 & 68.4& 71.4&25.5&13.1&\textbf{5.2}&11.6&10.8&\textbf{4.9}\\
    9 &131.7&137.3&50.7&22.3&\textbf{8.6}&21.1&19.9&\textbf{8.7}\\
    \midrule
    \multicolumn{9}{c}{Intel MKL 10.3.11.339 (intel64)} \\
    \midrule
    3 & 0.1& 0.1& 0.1& 0.1&0.1& 0.1& 0.1&\textbf{0.1}\\
    4 & 0.4& 0.3& 0.2& 0.2&0.2& 0.2& 0.2&\textbf{0.1}\\
    5 & 2.9& 1.3& 0.6& 0.6&0.4& 0.6& 0.6&\textbf{0.3}\\
    6 & 8.2& 4.2& 1.4& 2.1&0.7& 1.5& 1.1&\textbf{0.5}\\
    7 &20.5&16.3& 3.4& 5.9&1.7& 3.7& 2.4&\textbf{1.0}\\
    8 &35.7&34.1& 7.2& 9.1&3.7& 7.1& 5.5&\textbf{1.7}\\
    9 &73.3&73.8&14.4&17.6&5.9&12.8&11.2&\textbf{2.8}\\
    \bottomrule
  \end{tabular}
  \caption{M2L timings for bbFMM (prolate sphere)}
  \label{tab:timebbfmm_e2}
\end{table}

\begin{table}[htbp]
  \centering
    \begin{tabular}{c|ccc|cc|ccc}
    \toprule
    $Acc$ & NA & NAsym & NAblk & SA & SArcmp & IA & IAsym & IAblk \\
    \midrule
    \multicolumn{9}{c}{libblas 1.2.20110419-2} \\
    \midrule
    3 &  6.3&  5.9&  4.8& 2.0&\textbf{ 2.0}& 3.2& 2.9& 3.5\\
    4 & 25.7& 25.1& 20.1& 4.3&\textbf{ 3.9}& 8.0& 8.1& 8.3\\
    5 &113.0& 89.5& 71.4& 7.4&\textbf{ 7.1}&19.8&19.4&19.7\\
    6 &    -&275.9&202.2&14.5&\textbf{11.6}&43.0&42.8&40.4\\
    7 &    -&    -&    -&   -&            -&   -&86.5&78.1\\
    \midrule
    \multicolumn{9}{c}{Intel MKL 10.3.11.339 (intel64)} \\
    \midrule
    3 &  4.8&  3.9& 2.3& 1.6&\textbf{ 1.7}& 3.2& 2.4& 2.2\\
    4 & 21.5& 17.1& 7.5& 3.3&\textbf{ 3.4}& 6.7& 5.9& 4.8\\
    5 &109.1& 61.5&22.8& 6.3&\textbf{ 6.3}&10.6&14.0&10.0\\
    6 &    -&162.2&58.8&11.6&\textbf{ 9.8}&35.1&29.7&17.6\\
    7 &    -&    -&   -&   -&            -&   -&60.8&30.8\\
    \bottomrule
  \end{tabular}
  \caption{M2L timings for dFMM (sphere); high-frequency leaf level} 
  \label{tab:timem2l_e0}
\end{table}
\begin{table}[htbp]
  \centering
    \begin{tabular}{c|ccc|cc|ccc}
    \toprule
    $Acc$ & NA & NAsym & NAblk & SA & SArcmp & IA & IAsym & IAblk \\
    \midrule
    \multicolumn{9}{c}{libblas 1.2.20110419-2} \\
    \midrule
    3 & 1.8&  1.8&  1.6& 0.6&\textbf{ 0.6}& 0.9&  0.9& 1.1 \\
    4 & 8.9&  8.6&  6.8& 1.3&\textbf{ 1.1}& 2.7&  2.6& 2.7 \\
    5 &31.4& 30.6& 24.8& 2.8&\textbf{ 2.3}& 6.9&  6.8& 6.9 \\
    6 &91.8& 95.9& 71.6& 5.5&\textbf{ 4.1}&15.9& 15.9&14.8 \\
    7 &   -&228.4&176.6& 9.6&\textbf{ 7.3}&32.1& 32.3&28.9 \\
    8 &   -&    -&    -&16.4&\textbf{11.2}&59.2& 59.6&52.7 \\
    9 &   -&    -&    -&25.5&\textbf{21.9}&   -&105.1&90.1 \\
    \midrule
    \multicolumn{9}{c}{Intel MKL 10.3.11.339 (intel64)} \\
    \midrule
    3 & 1.4&  1.0& 0.7& 0.4&\textbf{ 0.4}& 0.7& 0.8& 0.6\\
    4 & 7.2&  5.0& 2.4& 0.9&\textbf{ 1.0}& 2.0& 1.8& 1.3\\
    5 &25.8& 20.7& 7.9& 2.2&\textbf{ 1.9}& 5.5& 4.5& 3.1\\
    6 &62.5& 56.7&20.4& 4.3&\textbf{ 3.4}&12.6&10.3& 5.9\\
    7 &   -&141.3&48.8& 7.8&\textbf{ 5.5}&24.7&21.8&11.0\\
    8 &   -&    -&   -&13.6&\textbf{ 9.5}&45.3&42.0&18.8\\
    9 &   -&    -&   -&21.0&\textbf{14.7}&   -&79.6&30.9\\
    \bottomrule
  \end{tabular}
  \caption{M2L timings for dFMM (oblate sphere); high-frequency leaf level} 
  \label{tab:timem2l_e1}
\end{table}
\begin{table}[htbp]
  \centering
    \begin{tabular}{c|ccc|cc|ccc}
    \toprule
    $Acc$ & NA & NAsym & NAblk & SA & SArcmp & IA & IAsym & IAblk \\
    \midrule
    \multicolumn{9}{c}{libblas 1.2.20110419-2} \\
    \midrule
    3 &  0.6 &  0.6 &  0.5 &  0.3 &\textbf{ 0.2}& 0.3 &  0.3 & 0.3 \\
    4 &  3.2 &  2.8 &  2.3 &  1.0 &\textbf{ 0.5}& 0.9 &  0.9 & 1.0 \\
    5 & 11.9 & 10.0 &  8.7 &  2.7 &\textbf{ 1.2}& 2.4 &  2.4 & 2.5 \\
    6 & 32.5 & 30.8 & 25.1 &  6.5 &\textbf{ 2.6}& 6.0 &  5.6 & 5.7 \\
    7 & 80.3 & 77.0 & 61.1 & 13.2 &\textbf{ 5.1}&12.1 & 11.8 &11.1 \\
    8 &    - &    - &    - & 24.8 &\textbf{ 9.1}&24.3 & 23.3 &21.2 \\
    9 &    - &    - &    - & 47.5 &\textbf{18.6}&43.5 & 43.5 &37.2 \\
    \midrule
    \multicolumn{9}{c}{Intel MKL 10.3.11.339 (intel64)} \\
    \midrule
    3 &  0.3 &  0.3 &  0.2 &  0.2 &\textbf{ 0.2}& 0.2 &  0.2 & \textbf{ 0.2} \\
    4 &  2.1 &  1.4 &  0.7 &  0.5 &\textbf{ 0.4}& 0.6 &  0.6 & \textbf{ 0.4} \\
    5 &  8.7 &  5.8 &  2.5 &  2.0 &\textbf{ 1.0}& 2.0 &  1.5 & \textbf{ 1.0} \\
    6 & 20.2 & 17.9 &  6.7 &  5.1 &\textbf{ 2.0}& 4.6 &  3.4 & \textbf{ 2.0} \\
    7 & 48.4 & 46.5 & 16.3 & 10.3 &\textbf{ 4.0}& 9.6 &  7.5 & \textbf{ 3.8} \\
    8 &    - &    - &    - & 16.6 &\textbf{ 7.2}&18.3 & 15.8 & \textbf{ 7.0} \\
    9 &    - &    - &    - & 31.6 &\textbf{14.8}&33.7 & 32.1 & \textbf{11.6} \\
    \bottomrule
  \end{tabular}
  \caption{M2L timings for dFMM (prolate sphere); low-frequency leaf level} 
  \label{tab:timem2l_e2}
\end{table}

\paragraph{Impact of using symmetries on the memory requirement}
Missing times in the tables indicate that the required memory exceeded the
available memory of $\unit[8]{GB}$. Clearly, the NA variant (no approximation,
no use of symmetries) requires the most memory. On the other hand, IAsym and
IAblk are the most memory friendly variants. Both require only memory for
storing the low-rank representations of M2L operators from $T_{\text{sym}}$
(see Tab.~\ref{tab:size_T_all} for a comparison of $T$ and $T_{\text{sym}}$).

\paragraph{Impact of the blocking scheme on runtime}
Bold numbers indicate the fastest variants. Two results attract our
attention. 1) If we look at the times in the Tabs.~\ref{tab:timebbfmm_e0}--\ref{tab:timem2l_e2} we notice that in four cases (bbFMM with the sphere, the
oblate sphere and the prolate sphere and dFMM with the prolate sphere) IAblk,
and in two cases (dFMM with the sphere and the oblate sphere) SArcmp is
fastest. To be more general, IAblk wins at levels having non-directional
expansions (all levels of bbFMM and low-frequency levels of dFMM) and SArcmp
at levels having directional expansions (high-frequency levels of dFMM). Why
is that? The reason follows from the cost studies in
Sec.~\ref{sec:recompression_ucb}. Let us take for example $Acc=5$. Recall, we
measure the computational cost based on the size of the approximate M2L
operators, ie., $\op{O}(k^2)$ for SA and $\op{O}(2k\ell^2)$ for IA. The
respective ranks $k$ are given in Tab.~\ref{tab:lowranks}. The ratio of these
costs for SA to IA is $0.46$ at the high-frequency level $6$ and $1.60$ at the
low-frequency level $7$. As a matter of fact, at high-frequency levels wins SA
and at low-frequency levels IA. Even savings of $47\%$ at the low-frequency
level due to the recompression of SArcmp (see Tab.~\ref{tab:flops_ucbe2}) are
not sufficient to outweigh the advantage of the blocking scheme of IAblk.
If there is no low-frequency level, such as for the sphere in
Tab.~\ref{tab:timem2l_e0} and the oblate sphere in Tab.~\ref{tab:timem2l_e1},
the SArcmp outperforms all other variants. For example, if we would repeat the
computations for the prolate sphere with an octree of depth $6$ (no
low-frequency level) the resulting timing patterns would follow those from the
sphere and the oblate sphere (the overall application time would increase too,
since the tree-depth is based on our choice of balancing near- and far-field,
i.e., the shortest overall application time).  2) Evident is the wide margin
in the speedup of variants that use blocking and those which do not. If we use
the MKL (as opposed to libblas) for NA, NAsym, SA, SArcmp, IA and IAsym we end
up with $1.5 - 2$ times faster application times. However, if we use the MKL
for NAblk and IAblk we achieve $3 - 4$ times faster times. Even though these
speedups are greatest with the MKL library, it highlights the benefits of the
blocking scheme presented in Alg.~\ref{alg:blockedm2l}.

\begin{figure}[h!]
  \centering
  \begin{tikzpicture}
    \begin{semilogyaxis}[
      ylabel = M2L application time $\unit{[s]}$,
      xlabel = Accuracy $Acc$,
      legend style={at={(1.02,1.0)}, font=\scriptsize, anchor=north west},
      ]
      \addplot[solid,color=green,mark=o] coordinates {
        (3,0.3)
        (4,1.2)
        (5,6.3)
        (6,14.1)
        (7,57.7)
        (8,117.0)
        (9,260.5)
      };
      \addplot[solid,color=green,mark=star] coordinates {
        (3,0.1)
        (4,0.3)
        (5,1.3)
        (6,4.2)
        (7,16.3)
        (8,34.1)
        (9,73.8)
      };
      \addplot[solid,color=blue,mark=o] coordinates {
        (3,0.2)
        (4,0.6)
        (5,2.0)
        (6,6.9)
        (7,19.3)
        (8,34.3)
        (9,60.4)
      };
      \addplot[solid,color=blue,mark=star] coordinates {
        (3,0.1)
        (4,0.2)
        (5,0.6)
        (6,2.1)
        (7,5.9)
        (8,9.1)
        (9,17.6)
      };
      \addplot[solid,color=red,mark=o] coordinates {
        (3,0.4)
        (4,0.7)
        (5,1.4)
        (6,2.6)
        (7,5.8)
        (8,12.0)
        (9,20.6)
      };
      \addplot[solid,color=red,mark=star] coordinates {
        (3,0.1)
        (4,0.2)
        (5,0.4)
        (6,0.7)
        (7,1.7)
        (8,3.7)
        (9,5.9)
      };
      \addplot[solid,color=brown,mark=o] coordinates {
        (3,0.4)
        (4,0.8)
        (5,2.0)
        (6,4.0)
        (7,8.2)
        (8,16.2)
        (9,32.9)
      };
      \addplot[solid,color=brown,mark=star] coordinates {
        (3,0.1)
        (4,0.2)
        (5,0.6)
        (6,1.1)
        (7,2.4)
        (8,5.5)
        (9,11.2)
      };
      \legend{NAsym (s),
        NAsym (ps), SA (s),
        SA (ps), SArcmp (s),
        SArcmp (ps), IAsym (s),
        IAsym (ps)}  
    \end{semilogyaxis}
  \end{tikzpicture}
  \caption{Running times versus accuracy for NAsym, SA, SArcmp and IAsym
    for bbFMM taken from the Tabs.~\ref{tab:timebbfmm_e0},
    \ref{tab:timebbfmm_e1} and \ref{tab:timebbfmm_e2}; (s) stands for sphere 
    and (ps) for prolate sphere}
  \label{fig:growth_appl}
\end{figure}
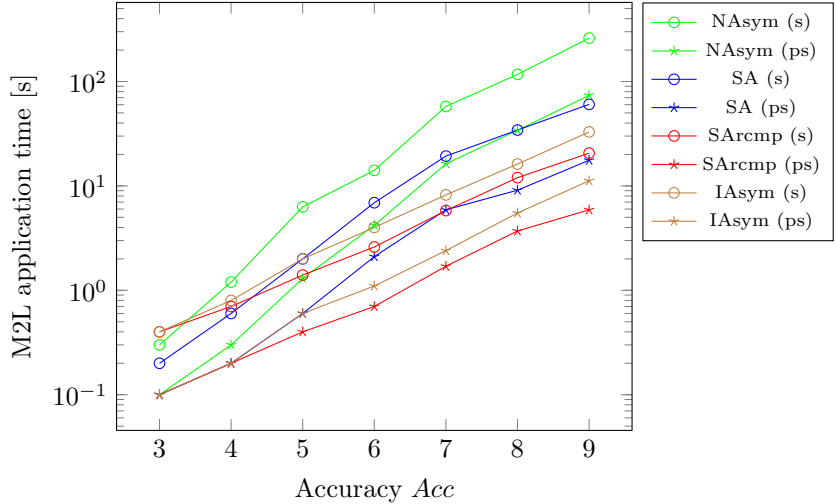

\paragraph{Varying growth of application times}
In Fig.~\ref{fig:growth_appl} we visualize the different growth rates of M2L
application times for the bbFMM with increasing accuracies $Acc$. We are
interested in the growth rates due to algorithmic changes. That is why we only
study those variants that do not use blocking. Since no approximation is
involved the times for NAsym grows the fastest. The times for SA grow slower
but still faster than those for IAsym. SArcmp features the slowest growth, it
is the optimal variant in terms of computational cost (see
Tab.~\ref{tab:flops_ucbe2}).


\section{Conclusion}
\label{sec:conclusion}

The fast multipole method based on Chebyshev interpolation, first presented in \citep{fong09a} for smooth kernel functions (bbFMM) and extended in \citep{Messner2011} to oscillatory ones (dFMM), is a convenient-to-use algorithm due to its easy implementation and its kernel function independence. In this paper we investigated algorithms to reduce the running time of the M2L operator. We proposed several optimizations and studied their respective strengths and weaknesses.

On one hand we proposed SArcmp, which uses an individual recompression of the
suboptimally approximated M2L operators obtained via SA (the variant presented
in \citep{fong09a}). We have shown that this new variant reduces the
computational cost noticeably. In some settings it even provides the fastest
M2L application times. On the other hand we also proposed a new set of
optimizations based on an individual low-rank approximation of the M2L
operators; we refer to them as IA variants. As opposed to SA they directly
lead to the optimal low-rank representation for each operator. The overall
number of flops is greater than for SArcmp (which is strictly a lower bound on
the number of flops). However, the advantage of the individual treatment of
the M2L operators is that we can exploit symmetries in their arrangement. This
means that all operators can be expressed as permutations of a subset. For
example, in the case of the bbFMM (in which the full interaction list has a
constant size), we need to approximate and store only $16$ instead of $316$
operators. The remaining ones can be expressed as permutations thereof. This
has a great impact on the precomputation time and the memory
requirement. Moreover it allows to express (again in the case of the bbFMM)
the at most $189$ matrix-vector products (applications of the M2L operators)
as at most $16$ {\bf matrix-matrix} products. We referred to this approach as
the IAblk variant. It can then take advantage of highly optimized
implementations of matrix-matrix operations (e.g., the MKL \citep{intel_mkl}).

Let us conclude by comparing SArcmp and IAblk, the two variants that have the
fastest running times. IAblk wins if we compare precomputation time, required
memory and runtime at levels having non-directional expansions (bbFMM and
low-frequency levels in dFMM). SArcmp wins only if we compare the runtime at
levels having directional expansions (high-frequency levels in dFMM). However,
in order to identify the optimal variant we have to distinguish two potential
uses of the FMM as a numerical scheme to perform fast matrix-vector
products. 1) If we are interested in the result of a single matrix-vector
product, a quick precomputation is essential. However, 2) if we are looking
for the iterative solution of an equation system (many matrix-vector
products), a fast running time of the M2L operator is crucial. Let us explain
this with an example. We take {\bf dFMM} (with MKL) with accuracy $Acc=5$ for
the sphere. IAblk wins if we are interested in the former use. The
precomputation takes $\unit[0.4]{s}$ versus $\unit[69.1]{s}$ (for SArcmp) and
the M2L application takes $\unit[10.0]{s}$ versus $\unit[6.3]{s}$, which sums
up to $\unit[10.4]{s}$ versus $\unit[75.4]{s}$. All other operators (P2P, P2M,
M2M, L2L and L2P) have nearly the same runtime in both cases, and their
runtimes are negligible compared to M2L. Looking at the latter use, SArcmp
starts being faster if the iterative solution requires more than $19$
matrix-vector products. For higher accuracies this threshold rises, e.g., for
$Acc=6$ it lies at $26$ matrix-vector products. In the case of {\bf bbFMM},
IAblk is always optimal.


\bibliographystyle{plainnat}
\bibliography{main}

\end{document}